\documentclass[12pt]{article}

\usepackage{amsmath,amssymb}
\usepackage{amsthm}
\usepackage{mathabx}
\usepackage{graphicx,tikz}
\usepackage{enumerate}
\usepackage{titlesec}
\usepackage{caption}
\usepackage{xcolor}
\usepackage{accents}

\newtheorem{Thm}{Theorem}[section]
\newtheorem*{Thm*}{Theorem}
\newtheorem{Prop}[Thm]{Proposition}
\newtheorem*{Prop*}{Proposition}
\newtheorem{Lem}[Thm]{Lemma}

\newtheorem{Cor}[Thm]{Corollary}
\newtheorem{Fact}[Thm]{Fact}

\newtheorem{Exa}[Thm]{Example}

\theoremstyle{remark}

\newtheorem{Rmq}[Thm]{Remark}

\theoremstyle{definition}
\newtheorem{Def}[Thm]{Definition}

\newtheorem*{Def*}{Definition}

\numberwithin{equation}{section}

\newcommand{\g}[1]{{\mathfrak #1}}
\newcommand{\m}[1]{\mathbb{ #1}}
\newcommand{\mc}[1]{\mathcal{ #1}}
\newcommand*{\dt}[1]{%
  \accentset{\mbox{\large\bfseries .}}{#1}}

     \def\ol{\overline}    
   
\def\al{\alpha}       \def\be{\beta}        
\def\de{\delta}       \def\eps{\varepsilon}  
\def\ka{\kappa}       \def\la{\lambda}      
                
\def\ph{\varphi}               
\def\om{\omega}              
             \def\Ph{\Phi}
                  \def\Om{\Omega}

\def\bt{\begin{Thm}}
	\def\et{\end{Thm}}
\def\br{\begin{Rmq}}
	\def\er{\end{Rmq}}

\def\bc{\begin{Cor}}
	\def\ec{\end{Cor}}
\def\bp{\begin{Prop}}
	\def\ep{\end{Prop}}
\def\bl{\begin{Lem}}
	\def\el{\end{Lem}}
\def\bd{\begin{Def}}
	\def\ed{\end{Def}}
\def\bq{\begin{quotation}}
	\def\eq{\end{quotation}}
\def\bfa{\begin{Fact}}
	\def\efa{\end{Fact}}
\def\bexa{\begin{Exa}}
	\def\eexa{\end{Exa}}

\def\ra{\rightarrow}
\def\vs{\vspace{1em}}

\begin{document}

\title{
Fourier transform in cyclic groups}
	\author{Yves Benoist
}
	\date{}
	
	\maketitle
	
\begin{abstract}
\noindent 
On a cyclic group of prime order, the non-trivial Dirichlet characters together with their Fourier transform
have constant modulus outside 0 and vanish at 0.  
Answering a question of H. Cohn, we construct new functions with these properties. The proof relies on Floer homology. We also apply this method to the  {\it biunimodular functions problem}.
\end{abstract}

	\renewcommand{\thefootnote}{\fnsymbol{footnote}} 
	\renewcommand{\thefootnote}{\arabic{footnote}}

\tableofcontents
	
\section{Introduction}
\label{secintexp}
	
The starting point of this article is the following 
problem asked by  Harvey Cohn in 1994 which is Problem 39 p.\,202 in \cite{MontgomeryCBMS}. The first aim of this article is to give an answer to this problem. The proof will rely not only  on tools coming from Cyclotomic Field Theory but also on tools coming from Floer Homology.

\subsection{The Harvey Cohn's problem}
\label{seccrival}

This problem is looking for a converse of a well-known property of  the Dirichlet characters which, by definition, are the characters of the multiplicative group $\m F_p^*$ of the prime field $\m F_p=\m Z/p\m Z$. We recall that a  character  $\chi: \m F_p^*\ra \m C^*$
is a map such that $\chi(k\ell)=\chi(k)\chi(\ell)$ for all
$k$, $\ell$ in $\m F_p^*$. We extend $\chi$ to $\m F_p$ 
by setting $\chi(0)=0$. When $\chi$ is non trivial, these functions $\chi$ satisfy Condition \eqref{eqnharveycohn} below.
More generally, we will also deal with cyclic groups $C_d:=\m Z/d\m Z$ of odd order $d$. 

\begin{Def*}
\label{defunifun}
A function $f:C_d\ra \m C$ such that $f(0)=0$ is said to be {\it unimodular} or {\it unimodular on $C_d\smallsetminus\{0\}$} if $|f(\ell)|=1$ for $\ell\neq 0$. 
\end{Def*} 

The function $f$ is said to be odd if $f(-k)=-f(k)$ for all $k$ in $C_d$.
\vs 
	
\noindent {\bf Problem (H. Cohn)} {\it Let 
$f:C_d\ra \m C$ be a function on the cyclic group $C_d$ vanishing at $0$ which is unimodular and satisfies
\begin{equation}
\label{eqnharveycohn}\textstyle
\sum\limits_{k\in C_d\smallsetminus\{ 0\}}f(k\!-\!\ell)\; f(k)^{-1}\; =\; -1
\;\;\;\;\mbox{for $\ell\neq 0$.}
\end{equation}
Assume $d\!=\!p$ is prime. Is the function  $f$  proportional to a Dirichlet character?}\vs

\noindent {\bf Geometric interpretation:} Condition \eqref{eqnharveycohn} means that the correlations between the unimodular function $f$ and its translates are constant equal to $\frac{-1}{d-1}
$. See also Lemma \ref{lemfoucohn} for a neat reinterpretation of this condition in terms of Fourier transform.

\begin{Def*}
A function $f$ on the cyclic group $C_d$
vanishing at $0$ and satisfying Condition \eqref{eqnharveycohn}
will be called a $\mc C$-function.
\end{Def*}

Before stating our results we gather a few comments. 
This problem is analogous to the problem of {\it biunimodular functions} that was raised by Per Enflo in the late 80's and quickly solved by Bj\"{o}rck and Saffari in \cite{BjorckSaffari}. 
See Appendix \ref{secbiufun} 
for a few words and new results on this analogous problem.

A counterexample to an extension of  Cohn's problem, where the finite field $\m F_p$ is replaced by a non-prime field $\m F_{p^k}$ with $k\geq 2$ 
and $p^k>4$, was found by Choi and Siu in 2000 in \cite{ChoiSiuJNT}.
Partial positive answers for prime fields were given 
by Borwein, Choi and Yazdani in 2001 in \cite[Cor. 1.2]{BorweinChoiYazdani}, by Kurlberg in 2002 when the values of $f$ are roots of unity in \cite{KurlbergJNT},   and by
Klurman and Mangerel in 2018 in \cite[Section 1.3]{KlurmanMangerel}
who refered to this problem as ``Cohn's conjecture''. 

When $d=p$ is prime, the set of $\mc C$-functions  such that $f(1)=1$ is finite. This was proven by Biro
in 1999 in \cite{BiroJNT}. 

\subsection{Main results} 

When $d$ is not prime, no unimodular $\mc C$-functions on 
$C_d$ were known to exist.
In the following proposition, we remark that they always exist.

\bp
\label{proexico0}
Let $d\geq 3$ be an odd integer.
Then there exist odd  unimodular $\mc C$-functions on the cyclic group $C_d$.
\ep

The first non trivial case is when $d=9$. 
We will see in Section  \ref{secexad09} that 
there are exactly $12$ odd unimodular $\mc C$-functions on $C_9$. 
 All of them are Galois conjugate.

When $d=p$ is prime, Proposition \ref{proexico0} is not useful since one already knows that the odd Dirichlet characters are unimodular $\mc C$-functions. However, in this case, we have a stronger existence result that  answers  Harvey Cohn's problem.

\bt 
\label{thmexicoh}
Let $p\geq 11$ be a prime number.
Then there exist odd uni\-modular $\mc C$-functions $f$ on $\m F_p^*$ 
which are not proportional to a Dirichlet character.
\et

One can check that, for $p=3$, $5$ and $7$, every unimodular  
$\mc C$-function on $\m F_p$ 
is  proportional to a Dirichlet character.

The first non trivial case of Theorem \ref{thmexicoh} is when $p=11$. We will see
in Section  \ref{secexad11} that  there are exactly 
$15$ odd unimodular $\mc C$-functions on $\m F_{11}$, all of them
except one  having multiplicity $1$. The exception is the Legendre character. Among these functions there are $5$ odd Dirichlet characters.
The remaining $10$ functions are Galois conjugate.

The functions $f$ on $C_d$ constructed both in 
Proposition \ref{proexico0} and in Theorem \ref{thmexicoh} are more delicate to understand than the Dirichlet characters. Indeed, even when $f(1)=1$, they may take their values in non-cyclotomic number fields.
\vs 

Forgetting about the unimodularity condition, we will be able to compute
the number of odd $\mc C$-functions.

\bt
\label{thmcouodd}
Let $p$ be an odd prime, then, counted with multiplicities,
the number of odd  $\mc C$-functions $f$ on $\m F_p$ with $f(1)=1$ is equal to the binomial coefficient $\binom{2k}{k}$ where $k:=(p-3)/2$.
\et

Note that, when $p\equiv 3$ mod $4$, the multiplicity for the Legendre character $f=\chi_0$
is equal to $2^{(p-1)/2}$. This will be proven in Proposition \ref{protrach0}.$a$.

\subsection{Biunimodular functions}
\label{secintbiu}

In the appendix, we will  apply  the same method to a popular 
analogous problem:  the construction of biunimodular functions
on $\m F_p$. 
Those are functions of constant modulus whose Fourier transform also has constant modulus.
We will construct in Theorem \ref{thmexibiu} {\it new} biunimodular functions on $\m F_p$. Here {\it new} means that the function is neither a gaussian function nor 
the translate of an $(\m F_p^*)^2$-invariant function.

\subsection{Organization}

Proposition \ref{proexico0} will be proven in Chapter \ref{secsymgeo}.
Both Proposition \ref{proexico0} and Theorem \ref{thmexicoh} will rely on
results in symplectic geometry dealing with intersection of a Lagrangian submanifold and its translate by a Hamiltonian diffeomorphism in a compact symplectic manifold. These results,
Theorems \ref{thmcho} and \ref{thmchobis}, are due 
to Biran, Entov and Polte\-rovich in \cite{BiranEntovPolterovich}
and to  Cheol-Hyun Cho in \cite{ChoClifford}. Their proof relies on Floer homology. 
\vs 

\noindent
Theorem \ref{thmexicoh} will be proven in Chapter \ref{secconequ}
for  non safe primes $p$ , i.e. when  $(p-1)/2$ is non prime.
In its proof, it will be crucial to deal with equivariant 
functions on $\m F_p$. This will be useful  to check that the known intersections of the corresponding Clifford torus with its image under the Fourier transform are transverse. This transversality will rely 
on a property of Jacobi sums given in Proposition \ref{projacsum}.
\vs 

\noindent 
Proposition \ref{projacsum} will be proven 
in Chapter \ref{seccycfie}. It will rely on a classical theorem in cyclotomic field theory: the  Stickelberger formula
for the Jacobi sums.
\vs 

\noindent
Theorem \ref{thmexicoh} will be proven in Chapter \ref{secequsaf}
for  safe primes $p$ , i.e. when  $(p-1)/2$ is prime.
The extra difficulty in this case is that we have to deal with 
non transverse tori and to estimate the multiplicity of the intersection at the Legendre character. 
We will then slightly perturb the Fourier transform so that the known intersections of the tori become transverse.
\vs 

\noindent
Theorem \ref{thmcouodd} will be proven in Chapter \ref{seccoucoh}.
There we will deal with more general equivariant functions, see Theorem \ref{thmcouequ}. In the proof, we adapt a strategy developed
by Haagerup in \cite{HaagerupFourier} for counting cyclic $d$-roots. 
The key ingredient is the little theorem of Chebotarev
(Lemma \ref{lemlitche})
as reformulated by Biro and Tao.
\vs 

\noindent
A few examples will be given in Chapter \ref{secexalis} that emphasize
the complexity of these new $\mc C$-functions.
\vs 

\noindent
Finally we will discuss  in Appendix \ref{secbiufun} the adaptation of these methods to the  analog problem of biunimodular functions. 
\vs 

I would like to thank A. Mangerel that pointed out to me the Harvey Cohn problem while I was lecturing on \cite{CSAGI} and \cite{CSAGII}.

\section{Using Symplectic Geometry}
\label{secsymgeo}

In this Chapter we prove Proposition \ref{proexico0}.

We begin in Section \ref{secunigro} by stating a theorem 
dealing with the unitary group, its action on the complex projective space and the Clifford tori. 

We interpret in Section \ref{secfoutra} the $\mc C$-functions 
via the Fourier transform. 

We prove in Section \ref{secexicoh} 
the existence of unimodular $\mc C$-functions on $C_d$.

\subsection{The unitary group}
\label{secunigro}

We recall that the complex projective space $\m C\m P^{n-1}$ 
is the set of lines of $\m C^{n}$ that we denote $p=[z_1,\ldots, z_n]$.
A Clifford torus is a compact  $(n\!-\!1)$-dimensional torus of the form
$\m T^{n-1}:=\{p=[z_1,\ldots, z_n]\in \m C\m P^{n-1}\mid |z_i|=1
\; \mbox{\rm for all}\; i \}$ in a unitary basis of $\m C^n$.
The unitary group $U=U(n)$ acts naturally on $\m C\m P^{n-1}$.

\bt 
\label{thmcho} {\bf (Biran, Entov, Polterovich, Cho)}
Let $\m C\m P^{n-1}$ be the complex projective space 
and $\m T^{n-1}\subset \m C\m P^{n-1}$ be the Clifford torus.
Then for all unitary transformation $u\in U$, 
one has $\m T^{n-1}\cap u\m T^{n-1}\neq\emptyset$.
\et

This theorem is due 
to Biran, Entov and Polterovich in \cite{BiranEntovPolterovich}
and to  Cheol-Hyun Cho in \cite{ChoClifford}. Both proofs rely on Floer homology. 
The key remark being that  $\m C\m P^{n-1}$  is a closed
symplectic manifold, that $\m T^{n-1}$ is a closed lagrangian submanifold and that the unitary transformation $u$ is a hamiltonian diffeomorphism of
$\m C\m P^{n-1}$. These four authors consider a closed Lagrangian submanifold $L$ in a closed symplectic manifold.
Under some extra assumption on $L$, for instance when $L$ is 
``monotone'',
they prove that $L$ cannot be displaced from itself by a hamiltonian diffeomorphism.
Therefore,  the Clifford torus $\m T^{n-1}$ in the projective space $\m C\m P^{n-1}$ cannot be displaced from itself by a unitary operator $u$ in $U(n)$.

Idel and Wolf reformulated in  \cite{idelwolf} Theorem \ref{thmcho} as
a decomposition theorem for the unitary group $U=U(n)$.

Let $p_0=\m Cv_0$ be the point on the Clifford torus $\m T^{n-1}$ where $v_0$ is the vector $v_ 0=(1,\ldots,1)$.
Let $V$  be  the stabilizer
$V:=\{ u\in U\mid u(v_0)=v_0\}$, and let $T\subset U$ be the maximal torus subgroup $D:=\{{\rm diag}(u_1,\ldots,u_n)\in U\}$.

\bc 
\label{corcho} {\bf (Idel, Wolf)}
One has the equality $U=DVD$.
\ec
This means that every unitary matrix $u$ can be decomposed as a product of three unitary matrices $u=d_1vd_2$ with both $d_i$ diagonal and with 
$\sum_j v_{ij}=1$ for all $i=1,\ldots, n$. Note that this decomposition is not {\it unique modulo the center of} $U$. See \cite{AnderssonBengtsson} for some examples.

The following theorem in \cite{ChoClifford} is more precise. 
We recall that an intersection point $p\in \m T^{n-1}\cap u\m T^{n-1}$
is transverse if the tangent spaces at $p$ intersect transversally, that is if $T_p\m T^{n-1}\cap T_p \, u\m T^{n-1}=\{0\}$.

\bt 
\label{thmchobis} {\bf (Cho)}
Let $\m C\m P^{n-1}$ be the complex projective space of $\m C^n$, let 
$\m T^{n-1}\subset \m C\m P^{n-1}$ be the Clifford torus,
and $u\in U(n)$ a unitary transformation.
If the intersection $\m T^{n-1}\cap u\m T^{n-1}$
is transverse,  
it contains at least $2^{n-1}$ points.
\et

\subsection{The finite Fourier transform}
\label{secfoutra}
	
	The Harvey Cohn problem can be formulated in terms of Fourier transforms.
	
	In the present paper we focus on the cyclic group $C_d:=\m Z/d\m Z$ where $d$ is an odd integer, which  will often be assumed to be prime $d=p$. 
	We will use the Fourier transform $f\mapsto Ff=\widehat{f}$ on
	$\m C^{C_d}$, with the sign convention given by,
	for all  $k\in C_d$,
\begin{equation}
\label{eqnfoufor}
\textstyle
\widehat{f}(k)= \frac{1}{\sqrt{d}}\sum\limits_{\ell\in C_d}e^{2i\pi k\ell/d}\;f(\ell).
\end{equation}
See for instance \cite{TerrasFourier}. 

A function $f:C_d\ra \m C$ is said to be 
{\it biunimodular} on $C_d\smallsetminus\{0\}$ if both 
$f$ and $\widehat{f}$ vanish at $0$ and are unimodular on $C_d\smallsetminus\{0\}$ as defined in  Section \ref{seccrival}.\vs 
	
	The following lemma is a straightforward remark that gives a clean interpretation of the $\mc C$-functions  in term of the Fourier transform.

\bl
\label{lemfoucohn}
$a)$ A function $f: C_d\ra \m C$ vanishing at $0$ is a $\mc C$-function if and only if there exists a function $g: C_d\ra \m C$ 
	vanishing at $0$ such that 
\begin{equation} 
\label{eqnfoucohn}
fg={\bf 1}_{C_d\smallsetminus\{0\}}
\;\;{\rm and}\;\;
\widehat{f}\,\widehat{\widecheck{g}}={\bf 1}_{C_d\smallsetminus\{0\}}.
\end{equation}
	
$b)$ In particular, a function $f: C_d\ra \m C$ is a unimodular $\mc C$-function if and only if $f$ is 
	biunimodular on $C_d\smallsetminus\{0\}$.
	\el
Here $\widecheck{g}$ is the function
given by $\widecheck{g}(x)=g(-x)$. Note that 
$\widehat{\widecheck{g}}=\widecheck{\widehat{g}}$.

\begin{proof}
$a)$ The function $\sqrt{d}\,\widehat{f}\,\widehat{\widecheck{g}}$ is the Fourier transform of the convolution 
$f*\widecheck{g}$ and the function 
$\sqrt{d}\,{\bf 1}_{C_d\smallsetminus\{0\}}$
is the Fourier transform of 
$d\de_{0}-{\bf 1}_{C_d}$.

Therefore, Equation \eqref{eqnfoucohn} 
is equivalent to 
$$
fg={\bf 1}_{C_d\smallsetminus\{0\}}
\;\;{\rm and}\;\;
f*\widecheck{g}=d\,\de_0-{\bf 1}_{C_d}
$$ 
which is nothing but Condition \eqref{eqnharveycohn}.

$b)$ Just apply Point $a)$ with the function $g=\ol{f}$.
\end{proof}

When $d=p$ is a prime number, every Dirichlet character $\chi$
on $\m F_p$	is a unimodular $\mc C$-function. 
Indeed, it satisfies the condition: $\chi\widehat{\chi}\equiv {\bf 1}_{F_p^*}$.
Moreover, when $\chi(-1)=-1$, this function $\chi$ 
	is odd.

\subsection{Existence of unimodular $\mc C$-functions}
\label{secexicoh}

\begin{proof}[Proof of  Proposition \ref{proexico0}]
Let $d= 2n+1$ be an odd integer and  $V_-$ be the vector space of odd functions on 
$C_d$. 
By using the basis $(E_j)_{1\leq j\leq n}$ of $V_-$ given by 
$E_j:=\de_j-\de_{-j}$, one identifies $V_-$ with $\m C^n$.
The Fourier transform 
$f\mapsto \widehat{f}$ 
is a unitary transformation of $V_-$ that we still denote by $F$.
The elements of the Clifford torus 
$\m T^{n-1}$ of 
$\m P(V_-)=\m C\m P^{n-1}$ are precisely the lines spanned by odd unimodular functions 
on $C_d\smallsetminus \{ 0\}$.
Theorem \ref{thmcho} tells us that $\m T^{n-1}\cap F(\m T^{n-1})\neq\emptyset$. 
By Lemma \ref{lemfoucohn}.$a$ this exactly means that there exists a unimodular odd $\mc C$-function.
\end{proof}

\section{Equivariant unimodular $\mc C$-functions}
\label{secconequ}

In this Chapter we begin the proof of Theorem \ref{thmexicoh}, relying on Proposition \ref{projacsum} that will be proven in 
Chapter \ref{seccycfie}.

We introduce in Section \ref{secequcoh} the space of
$(H,c)$-equivariant functions.

We give a criterion in Section \ref{sectracli} for the transversality of 
a Clifford torus and its image by Fourier transform 
at the points given by Dirichlet characters.

We explain in Section \ref{secconuni} the examples where this transversality is not satisfied, and, by avoiding them
we deduce the existence of new 
unimodular $\mc C$-functions when the prime $p$ is a safe prime.

\subsection{The Jacobi sums}
\label{secjacsum}

Let $p\geq 3$ be a prime number. We recall that a Dirichlet character
$\chi$ on $\m F_p$ is said to be {\it principal} or {\it trivial} if $\chi\equiv {\bf 1}$ on
$\m F_p^*$, and that the order of $\chi$ is the smallest positive integer $d_\chi$ such that $\chi^{d_\chi}$ is trivial.  The Dirichlet character
$\chi_0$ of order $2$ 
is called the Legendre character, since it is given by the Legendre symbol. Let $\chi_1$, $\chi_2$ be two Dirichlet characters.
The corresponding Jacobi sum is defined as 
\begin{equation} 
\label{eqnjacsum}
J(\chi_1,\chi_2)=\sum_{x\in \m F_p}\chi_1(x)\,\chi_2(1-x).
\end{equation}
These algebraic numbers 
live in the cyclotomic field $K=\m Q(\zeta_{p-1})$ 
spanned by the $(p\!-\! 1)^{\rm th}$-root of unity
$\zeta_{p-1}=e^{\frac{2i\pi}{p-1}}$

Recall that, one has the equalities 
\begin{equation*}
	\label{eqnj12j21}
	J(\chi_1,\chi_2)=J(\chi_2,\chi_1)=\chi_2(-1)J(\ol{\chi}_1\ol{\chi}_2,\chi_2)
\end{equation*}
When $\chi_1\chi_2$ is not principal, one has the equality
\begin{equation}
	\label{eqnjacgau}
	J(\chi_1,\chi_2)=\tfrac{G(\chi_1)G(\chi_2)}{G(\chi_1\chi_2)}.
\end{equation}
where $G(\chi)$ is the Gauss sum 
\begin{equation}
	\label{eqngausum}
	G(\chi)=\textstyle\sum\limits_{x\in \m F_p}\chi(x)\zeta_{p}^{x}.
\end{equation}
Therefore, when $\chi_1$, $\chi_2$ and $\chi_1\chi_2$ are not principal, one has 
\begin{eqnarray} 	
	\label{eqnj12srp}
	|J(\chi_1,\chi_2)|=\sqrt{p}.
\end{eqnarray}
These quantities are useful: the Jacobi sum occurs in the convolution product
$$
\chi_1 *\chi_2
\;=\; J(\chi_1,\chi_2)\, \chi_1\chi_2,
$$
and, for $\chi$ non trivial, the Gauss sum occurs in the Fourier transform
\begin{equation} 
	\label{eqnchieig}
	\widehat\chi=\tfrac{G(\chi)}{\sqrt{p}}\, \ol{\chi}.
\end{equation}
When $\chi={\bf 1}_{\m F_p^*}$ is trivial, one has an extra term:\;
$	\widehat{{\bf 1}_{\m F_p^*}}=\tfrac{-1}{\sqrt{p}}\, {\bf 1}_{\m F_p^*}+ \tfrac{p-1}{\sqrt{p}}\de_0\,.$

We will need the following

\bp 
\label{projacsum} 
Let $p\geq 3$ be a prime number and 
$\chi_1$, $\chi_2$ be two non principal Dirichlet characters. Then one has the equivalence:
\begin{equation}
\label{eqnratjac}
\mbox{the ratio}\;\;	
	R(\chi_1,\chi_2):=
	\frac{J(\ol{\chi}_1,\chi_2)}{J(\chi_1,\chi_2)}
\;\;\mbox{is a root of unity}
\end{equation}
if and only if one of the following seven cases is satisfied:\\
$(a)$ $\chi_1$ has order $2$.\\
$(b)$ $p\equiv 1$ {\rm mod} $6$, \; $\chi_1$ has order $3$ and $\chi_2$ has order $6$.\\
$(c)$ $p\equiv 1$ {\rm mod} $10$,\;  $\chi_1$ has order $5$ and $\chi_2$ has order $10$.\\
$(d)$ $p\equiv 1$ {\rm mod} $12$,\; $\chi_1$ has order $4$ and $\chi_2$ has order $6$.\\
$(e)$ $p\equiv 1$ {\rm mod} $12$,\; $\chi_1$ has order $3$ and $\chi_2$ has order $4$.\\
$(f)$ $p\equiv 1$ {\rm mod} $30$,\; $\chi_1$ has order $5$ and $\chi_2$ has order $6$.\\
$(g)$ $p\equiv 1$ {\rm mod} $30$,\; $\chi_1$ has order $3$ and $\chi_2$ has order $10$.
\ep

Proposition \ref{projacsum} is an output of Stickelberger's theorem which gives the prime factorization in $\m Z[\zeta_{p-1}]$
of the Jacobi sums. We will give a detailed proof in Chapter \ref{seccycfie}.
\vs

We will see in Proposition \ref{protracli} that the meaningful case is when the ratio \eqref{eqnratjac} is
equal to $\chi_2(-1)$. We will also see that this special value of the ratio 
does not play a role in the proof and that it can occur in each of these seven cases.

\subsection{Equivariant $\mc C$-functions}
\label{secequcoh}

Before going on let us explain why we need to deal with equivariant  functions. We want to use Theorem \ref{thmchobis} 
to prove the existence of $\mc C$-functions. The Fourier transform 
$F$ is 
a unitary operator on the space $\m C^{\m F_p}$  of functions $f$ on $\m F_p$,
but it does not preserve the subspace of functions vanishing at $0$. 
We would like to find a vector subspace $V$ of this subspace 
such that both $V$ and $F(V)$ have a unitary basis for which the sup norm is proportional to the sup norm on $\m C^{\m F_p}$. 
The space $V_{H,c}$ of equivariant  functions on $\m F_p$ will  play this role. 
\vs 

Let $p\geq 3$ be a prime number, let $C_p=\m F_p\simeq \m Z/p\m Z$ be the cyclic additive group, and let 
$G:=\m F_p^*\simeq \m Z/(p\!-\!1)\m Z$ be the cyclic multiplicative group. Let $H\subsetneq
G$ be a proper subgroup of index $n\geq 2$, 
and $c:H\ra \m C^*$ be a non-trivial character. 
We introduce the vector space $V=V_{H,c}$ of $(H,c)$-equivariant functions $f$ on $\m F_p$.
\begin{equation}
	\label{eqnequfun}
	V_{H,c}:=\{ f:\m F_p\ra \m C\mid f(hx)=c(h)f(x)\;
	\mbox{\rm for all}\; h\in H,\; x\in \m F_p\}.
\end{equation}
Note that all $(H,c)$-equivariant functions $f$  vanish at $0$: one has $f(0)=0$. 
We denote by $g_1=1,g_2,\ldots ,g_n$ a family of representatives in $G$ of the classes $gH$, and by $(f_i)_{1\leq i\leq n}$ the family of elements of $V$ defined by $f_i(g_j)=\de_{i,j}$ for all 
$i$, $j$. This family is a basis of $V$. More precisely any function $f$ in $V$ can be written in a unique way as 
$f=\sum_{1\leq i\leq n} \la_if_i$ with $\la_i\in \m C$ and one has 
the equality of the sup norms 
$$
\sup_{x\in \m F_p} |f(x)|= \sup_{1\leq i\leq n}
|\la_i|.
$$

\subsection{Transversality of Clifford tori}
\label{sectracli}

Let 
$
\m P(V_{H,c})\simeq \m C\m P^{n-1}
$ 
be the projective space  of $V_{H,c}$ and $T_{H,c}$ be the Clifford torus
\begin{equation}
\label{eqnclitor}
T_{H,c}:=\{[f]\in \m P(V_{H,c})\mid \;|f(x)|=|f(1)|\; \mbox{\rm for all $x\in \m F_p^*$}\}\simeq \m T^{n-1}.
\end{equation}
The vector  space $V_{H,c}$ contains exactly $n$ Dirichlet characters $\chi$. 
The points $[\chi]$ belong to $T_{H,c}$. 
Since the Fourier transform $F:f\mapsto \widehat{f}$ is a unitary transformation 
that sends  $V_{H,c}$ onto $V_{H,\ol{c}}$ where $\ol{c}$ is the conjugate character,
we will have to deal with two hermitian vector spaces instead of one. 
This will not be an issue when applying Cho's theorem.
According to \eqref{eqnchieig}, the points $[\chi]\in \m P(V_{H,c})$ 
belong to  $T_{H,c}\cap F^{-1}T_{H,\ol{c}}$.

\bp 
\label{protracli}
Let $\chi$ be a character of $G:=\m F_p^*$ that extends the non trivial character $c$ of a subgroup $H\subset G$. 
The intersection $T_{H,c}\cap F^{-1}T_{H,\ol{c}}$ is transverse at $[\chi]$ if and only if, for all non trivial character $\psi$ of $G$ which is trivial on $H$, one has 
\begin{equation} 
	\label{eqnjchjps}
	J(\chi,\psi)\neq \psi(-1)\,J(\ol{\chi},\psi).
\end{equation}
\ep

We will see  in Proposition \ref{procouni}  nice conditions on $(H,c)$ that ensure 
\eqref{eqnjchjps}. 

\br 
Note that the factor $\psi(-1)$ is nothing but a sign 
which is given by the parity 
of $\psi$. When  $H$ has even order, the element $-1$ belongs to $H$ and this sign is equal to $+1$. 
\er

\begin{proof}[Proof of Proposition \ref{protracli}]
	Let $d=p\! -\! 1$, and write $d=n \,d_H$ where $n$ is the index of $H$
	and $d_H$ the order of $H$.
	\vs 
	
	$\star$ {\bf We first describe the tangent space to the projective space.}
	We fix $\chi$ as in the proposition and we set
	\begin{eqnarray*}
		B=B_{H,c}&:=&\{ \chi' \;\mbox{\rm character of $G$}\mid \chi'\in V_{H,c} \},\\
		B_o:=B_{H,o}&:=&\{ \psi \;\mbox{\rm character of $G$}
		\; \mbox{\rm trivial on $H$}\}.
	\end{eqnarray*}
	The set $B$ is a  basis of $V_{H,c}$. Similarly, the set $B_o$ is a basis of the space $V_{H,o}$ of $H$-invariant functions on $\m F_p$ 
	that vanish at $0$. Moreover every element $\chi'$ of $B$ can be written in a unique way as $\chi'=\chi\psi$ with $\psi$ in $B_o$. 
	This basis $B$ is orthogonal: one has, for $\chi'$, $\chi''$ in $B$,
	$$
	\langle\chi',\chi''\rangle_{\ell^2(\m F_p)} =
	\left\{
	\begin{array}{ll}
		p\!-\! 1  & \mbox{if } \chi'=\chi'' \\
		0 & \mbox{otherwise.}
	\end{array}
	\right.
	$$
	We set $B'_o:=B_o\smallsetminus\{\psi_0\}$ where $\psi_0={\bf 1}_G\in B_o$ is the trivial character of $G$.
	It will be convenient to use the following coordinates system 
	${\bf a}=({\bf a}_\psi)_{\psi\in B'_o}$ of 
	$\m P(V)$ in the neighborhood of $[\chi]$ 
	where the coordinates ${\bf a}_\psi$ are complex numbers. It is given by 
\begin{equation}
\label{eqncooafa}
	{\bf a}\mapsto [f_{\bf a}]
	\;\;{\rm where}\;\; 
	f_{\bf a}=
	\textstyle\left({\bf 1}_G+\sum_{\psi\in B'_o}{\bf a}_\psi\psi\right)\chi\, .
\end{equation}
These coordinates ${\bf a}=({\bf a}_{\psi})\in \m C^{B'_o}$ are also a linear coordinate system
for the tangent space of $\m P(V)$ at the point $[\chi]$, thanks to the formula
\begin{equation*}
		\label{eqntanpv}
	{\bf a}\mapsto v_{\bf a}:=\frac{d}{d\eps}[f_{\eps{\bf a}}]|_{\eps=0}
	\;\in\; T_{[\chi]}\m P(V_{H,c}) .
\end{equation*}
	
	$\star$ 
	{\bf We describe the tangent space to the Clifford torus $T_{H,c}$}.
	The real linear equations defining the tangent space 
	of $T_{H,c}$ at the point $[\chi]$ are
	\begin{equation*}
		\label{eqntanthc1}
		\frac{d}{d\eps}\left.\left(| f_{\eps{\bf a}}(x)|^2-| f_{\eps{\bf a}}(1)|^2\right) \right|_{\eps=0} \;=\; 0\; , \;\;\;
		\mbox{\rm for all $x$ in $G$}.
	\end{equation*}
	Since $\psi(1)=1$ for all $\psi$ in $B_o$, using \eqref{eqncooafa}, this can be rewritten as
	\begin{equation}
		\label{eqntanthc2}
		\textstyle
		{\rm Re}(\sum_{\psi\in B'_0}{\bf a}_\psi (\psi -\psi_0)) \;=\; 0.
	\end{equation}
Since the basis $B_o$ is invariant by complex conjugation, Condition \eqref{eqntanthc2} can be rewritten as 
	\begin{equation*}
		\label{eqntanthc3}
		\textstyle
		\sum_{\psi\in B'_0}(\ol{{\bf a}_{\ol \psi}}+{\bf a}_\psi) \,
		(\psi -\psi_0) \;=\; 0.
	\end{equation*}
	By the linear independance of the characters of $G/H$, this gives
	\begin{equation}
		\label{eqntanthc4}
		T_{[\chi]}T_{H,c} \;\simeq\;
		\{( {\bf a}_{\psi})\in \m C^{B'_o}\mid 
		\ol{{\bf a}_{\ol \psi}}=-{\bf a}_\psi
		\; \;\mbox{\rm for all $\psi\in B'_o$}\}.
	\end{equation}
	
	$\star$ {\bf We describe the tangent space to the torus $F^{-1} T_{H,\ol{c}}$}.
	Using \eqref{eqnchieig}, one computes in our coordinate system
	$$
	\widehat{f}_{\bf a}=
	\textstyle
	\frac{G(\chi)}{\sqrt{p}}\left({\bf 1}_G+\sum_{\psi\in B'_o}\al_\psi {\bf a}_\psi\ol{\psi}\right)\ol{\chi}
	$$
	where we have, using also \eqref{eqnjacgau},
	\begin{equation}
		\label{eqnalpsi}
		\al_\psi
		\;:=\;
		\frac{G(\chi\psi)}{G(\chi)}
		\;=\;
		\frac{G(\psi)}{J(\chi,\psi)}\; .
	\end{equation} 
One deduces from \eqref{eqntanthc4} the equality
	\begin{equation}
		\label{eqntanthc5}
		T_{[\chi]}F^{-1}T_{H,\ol{c}} \;\simeq\;
		\{( {\bf a}_{\psi})\in \m C^{B'_o}\mid 
		\ol{\al_{\ol \psi}}\,\ol{{\bf a}_{\ol \psi}}=-\al_\psi\,{\bf a}_\psi
		\; \mbox{\rm for all $\psi\in B'_o$}\}.
	\end{equation}
		
	$\star$ {\bf We give the transversality criterion for the tangent spaces}.
	One easily computes
	\begin{equation}
		\label{eqnalpsb}
		\ol{\al_{\ol{\psi}}}
		\;=\;
		\frac{\ol{G(\ol{\psi})}}{\ol{J(\chi,\ol{\psi})}}\; 
		\;=\;
		\frac{\psi(-1)\,G(\psi)}{J(\ol{\chi},\psi)}\; .
	\end{equation} 
	Comparing \eqref{eqntanthc4} and \eqref{eqntanthc5}, and using 
	the values \eqref{eqnalpsi} and \eqref{eqnalpsb} for 
	$\al_\psi$ and $\ol{\al_{\ol{\psi}}}$, one gets the equivalences:
	\begin{eqnarray*}
		\label{eqncontra}
		T_{[\chi]}T_{H,c}\!\cap\! 	T_{[\chi]}F^{-1}T_{H,\ol{c}} =\{0\}
		&\Longleftrightarrow&
		\ol{\al_{\ol \psi}}\neq\al_\psi
		\;\; \mbox{\rm for all $\psi\in B'_o$}\\
		&\Longleftrightarrow&
		J(\chi,\psi)\!\neq\! \psi(-1)\,J(\ol{\chi},\psi)
		\;\; \mbox{\rm for all $\psi\in B'_o$}.
	\end{eqnarray*}
	This ends the proof of Proposition \ref{protracli}.
\end{proof}

\subsection{Unimodular $\mc C$-functions for non safe prime}
\label{secconuni}

In this section we explain how to check that
the intersection $T_{H,c}\cap F^{-1}T_{H,\ol{c}}$ is not transverse at 
a Dirichlet character $[\chi]$. We deduce from that the 
existence of new unimodular $\mc C$-functions for a non safe prime $p$.
\vs 

Let $p\geq 3$ be  prime. 
Let $H\subsetneq
\m F_p^*$ be a proper subgroup of index $n$
and $c:H\ra \m C^*$ be a non-trivial character
of order $d_c$. 
Let $V=V_{H,c}$ be the $n$-dimensional vector space of $(H,c)$-equivariant functions $f$ on $\m F_p$ as in
\eqref{eqnequfun}.

\bp
\label{procouni}
Assume that we are not in one of the  four cases\\
$(i)$  $d_c=2$ and $n$ is odd,\\
$(ii)$  $d_c=2$ and $n\equiv 6$ mod $12$,\\
$(iii)$ $d_c=3$ and $n$ is a multiple of $4$ or $10$ and $n$ is coprime to $3$,\\
$(iv)$ $d_c=5$ and $n$ is a multiple of $6$ and $n$ is coprime to $5$.\\
$a)$ Then the  intersection $T_{H,c}\cap F^{-1}T_{H,\ol{c}}$ is  transverse at 
$[\chi]$, for all Dirichlet character $\chi$
extending $c$.\\
$b)$ When $n\geq 3$, there exists in $V_{H,c}$ a unimodular $\mc C$-function which is not proportional to a Dirichlet character.
\ep

Here is a comment on Case $(i)$, and more precisely when the order of $\chi$ is $2$.
Recall that there is only one  Dirichlet character $\chi_0$ of order $2$. It is given by the Legendre symbol: 
$\chi_0(x) =(\!\frac{x}{p}\!)$ for all $x$ in $\m F_p$. 
It belongs to $V_{H,c}$ if and only if $d_c=2$ 
and $n$ is odd. In this case, the intersection $T_{H,c}\cap F^{-1}T_{H,\ol{c}}$ is not transverse at $[\chi_0]$. 
In this case it will be more delicate to apply Theorem \ref{thmchobis}.
This will be done in Chapter \ref{secequsaf}. 

\begin{proof}[Proof of Proposition \ref{procouni}]
$a)$ 
This follows from Propositions \ref{projacsum} and \ref{protracli}.
 We look at the lists of exceptions 
$(d_\chi,d_\psi)$ in  Proposition \ref{projacsum}
remembering that $d_\chi$ is a multiple of $d_c$ and that 
$d_\psi$ divides $n$. 
We also use the fact that,
when the index $n$ is not coprime to $d_c$, 
the character $c$ of $H$ does not admit an extension $\chi$ 
of order $d_c$.
In particular,  using the labeling of Proposition \ref{projacsum},
Case $(a)$ implies $(i)$,  Case $(b)$ cannot occur, 
Case $(c)$ also cannot occur, Case $(d)$ implies $(ii)$. 
Case $(e)$ implies $(iii)$, Case $(f)$ implies $(iv)$ and 
Case $(g)$ also implies $(iii)$.

$b)$ 	
Assume, by contradiction that the intersection $T_{H,c}\cap F^{-1}T_{H,\ol{c}}$ contains only Dirichlet characters $[\chi]$.
Then, by Point $a)$, this intersection is transverse. Therefore 
Theorem \ref{thmchobis} predicts the existence of 
at least $2^{n-1}$ intersection points.
Since $n\geq 3$, one has $2^{n-1}>n$. Since the number of Dirichlet characters in $V_{H,c}$ is $n$, there must exist another intersection point. This is the contradiction we are looking for.
\end{proof}
 
A prime $p$ for which $(p-1)/2$ is prime is called 
a safe prime. For instance $p=7$, $11$, $23$, $47$,  ...
Conjecturally, there are infinitely many safe primes, but,  
by Dirichlet Theorem, most of the primes are non safe.

\bc
\label{corthcohn}
Let $p\geq 11$ be a non safe prime.\\
$a)$ There exist a subgroup $H\subset \m F_p^*$ of index $n\!\geq\! 3$ containing $-1$
and an odd character $c$ of $H$ such that, for all Dirichlet character $\chi$ in $V_{H,c}$, 
the intersection $T_{H,c}\cap F^{-1}T_{H,\ol{c}}$ 
is transverse at $[\chi]$.\\
$b)$ For such a pair $(H,c)$ the space 
$V_{H,c}$  contains an odd  uni\-modular $\mc C$-function $f$
which is not proportional to a Dirichlet character.
\ec

\begin{proof}[Proof of Corollary \ref{corthcohn}]
$a)$ We recall that an odd character $c$ of $H$ is a character such that $c(-1)=-1$. For such a character, all the functions $f$ in $V_{H,c}$ are odd. We set $d=p\! -\! 1$ and distinguish two cases.

{\bf First case:} $d$ has an odd prime factor $\ell$.\\
We set $r:=d/\ell$. We choose $H$ to be the subgroup of $\m F_p^*$ of prime index $n:=\ell\geq 3$ and $c$ to be a character of $H$ of order $d_c:=r\geq 3$. 
Since these values do not occur in the exceptions of Proposition 
\ref{procouni} 
the intersection $T_{H,c}\cap F^{-1}T_{H,\ol{c}}$ is transverse at $[\chi]$.  This group $H$ contains $-1$ and one has $c(-1)=-1$.

{\bf Second case:} $d$ is a power of  $2$.\\
We set $r:=d/4$. It is also a power of $2$.
We choose $H$ to be the subgroup of $\m F_p^*$ of index $n:=4$ and $c$ to be the character of $H$ of order $d_c:=r$. 

Again, since these values do not occur in the exceptions of Proposition 
\ref{procouni} 
the intersection $T_{H,c}\cap F^{-1}T_{H,\ol{c}}$ is transverse at $[\chi]$.  This group $H$ contains $-1$ and one has $c(-1)=-1$.

$b)$ By the proof of Proposition \ref{procouni}, the space 
$V_{H,c}$  contains a uni\-modular $\mc C$-function $f$ on $\m F_p^*$ which is not proportional to a Dirichlet character.
\end{proof}

\br 
The conclusion of Proposition \ref{procouni}.$a$, i.e. the transversality of the intersection, is not true in  the exceptional cases $(i)-(iv)$ 
for infinitely many values of $p$. 
\er

\section{Using cyclotomic field theory}
\label{seccycfie}

The aim of this Chapter is to prove Proposition \ref{projacsum}.
I give a full proof since I did not find a reference for it. 
The key ingredient is the Kummer-Stickelberger formula which is a formula that give, both for Gauss sums and for Jacobi sums, an explicit
factorization as a product of prime ideals.

In Section \ref{secactzdz}, we begin  by a combinatorial lemma which parametrizes the orbits of $(\m Z/d\m Z)^*$ acting on $(\m Z/d\m Z)^2$.
The proof involves many cases.

In Section \ref{secstifor} we recall the short and 
elementary calculation which implies the Stickelberger's theorem.

In Section \ref{secdisjac}, we use this combinatorial lemma together with the Stickelberger theorem to distinguish various Jacobi sums.

In Section \ref{secexajac}, we point out surprising equalities  between Jacobi sums that prevent the corresponding Clifford tori 
to intersect transversally.

\subsection{Action on $(\m Z/d\m Z)^2$}
\label{secactzdz}

We will need the following combinatorial lemma.
Let $d\geq 2$ be an integer. 
We introduce an equivalence relation 
$(j,k)\sim (j',k') $ on $\m Z^2$ given by
\begin{equation}
\label{eqnequrel}
\begin{array}{c}
\mbox{\rm there exists an integer $x$ coprime to $d$ such that} \\
j' = xj\; {\rm mod} \; d\;\;\; {\rm and}\;\;\; k' = xk\; {\rm mod} \; d.
\end{array}
\end{equation}
We will see that these equivalence classes are related to the orbits of the Galois group acting on the Jacobi sums. 
The following lemma gives a nice representative in each equivalence class.

\bl 
\label{lemcomzdz}
Let $d\geq 2$ be an integer. 
For every positive integers $j,k<d$, there exist  positive integers $j',k'<d$ such that 
\begin{equation}
\label{eqnjkdjkd}
\begin{array}{l}
(j',k')\sim (j,k)\;\; {\rm or}\;\; (j',k')\sim(-j,k) \, ,\\
{\rm and}\;\;\;\;\; j'\leq k'\leq d -j',
\end{array}
\end{equation}
except in one of the following seven cases, where $m\geq 1$ is an integer:\\
$(a)$ $d=2\, m$,\; $j=m$\; and\;  $k\neq m$.\\
$(b)$ $d=6\,m$,\; 
$j\equiv \pm 2\,m$ {\rm mod} $d$\; and\; 
$k\equiv \pm m$ {\rm mod} $d$.\\
$(c)$ $d=10\,m$,
$j\equiv \pm 2\,m$ {\rm mod} $d$\; and\; 
$k\equiv \pm m$ {\rm mod} $d$.\\
$(d)$ $d=12\,m$,
$j\equiv \pm 3\,m$ {\rm mod} $d$\; and\; 
$k\equiv \pm 2\,m$ {\rm mod} $d$.\\
$(e)$ $d=12\,m$,
$j\equiv \pm 4\,m$ {\rm mod} $d$\; and\; 
$k\equiv \pm 3\,m$ {\rm mod} $d$.\\
$(f)$ $d=30\,m$,
$j\equiv \pm 6\,m$ or $\pm 12\, m$ {\rm mod} $d$\;\; and\;\; 
$k\equiv \pm 5\,m$ {\rm mod} $d$.\\
$(g)$ $d=30\,m$,
$j\equiv \pm 10\,m$\;\; and\;\; 
$k\equiv \pm 3\,m$\;
or $\pm 9\, m$ {\rm mod} $d$.
\el 

\br 
Note that in each of these seven cases $(a)$ to $(g)$, one can not find a pair $(j',k')$ satisfying \eqref{eqnjkdjkd}.
\er

For the proof we will need
the Jacobsthal function $g(n)$ which is  the largest gap in the sequence of  integers coprime to $n$. 
Among $g(n)$ consecutive integers, there is always one which is coprime to $n$.
This gap $g(n)$ grows very slowly with $n$. 
It can be bounded by a function that depend 
only on the number $\omega(n)$ of prime factors of $n$. 
Here, we will only need an upper bound on $g(n)$
which is due to Kanold.

\bl 
\label{lemkanold} {\bf (\cite{Kanold})}
For all integers $n\geq 2$, one has $g(n)\leq 2^{\,\omega(n)}$.\\
In particular one has $g(n)\leq n/3$ as soon as $n>10$.
\el

\begin{proof}[Proof of Lemma \ref{lemcomzdz}]
We assume that $j$, $k <d$ are positive integers 
for which one can not find integers $j',k'$ satisfying 
\eqref{eqnjkdjkd}, and we want to prove that we are in one of the 
seven cases $(a),\ldots,(g)$. 
\vs 

{\bf First step: We can assume that $d =j r$ with $r$ integer $r\geq 3$}.\\
Indeed there exists an integer $x$ coprime to $d$ such that $xj=j'$ mod $d$ and $j'$ is a positive divisor of $d$ so that $d=j' r$ with $r\geq 2$. If $r=2$, this is Case $(a)$. 
\vs 

{\bf Second step: We can assume that $j$ is coprime to $k$.}\\
Indeed if this is not the case we argue by induction.
We introduce the integer $m:={\rm gcd}(j,k)$ and set $j_0:=j/m$, $k_0:=k/m$
and $d_0:=d/m$. We find an equivalent pair $(j'_0,k'_0)$ satisfying 
\eqref{eqnjkdjkd} with the integer $d_0$. Then the pair 
$(j',k'):=(mj'_0,mk'_0)$ satisfies \eqref{eqnjkdjkd} with 
the integer $d$.
\vs 

{\bf Third step: We can assume that $k < j$.}\\
Indeed since \eqref{eqnjkdjkd} is not satisfied 
we have either $k<j$ or $k>d-j$. In the second case,
we replace the pair $(j,k)$ by $(j',k'):=(j,d-k)\sim (-j,k)$.
\vs 

{\bf Fourth step: We can assume that $j$ is coprime to $r$.}\\
Indeed let $s:=gcd(r,j)$ and $t:=lcm(r,j)=d/s$. Assume by contradiction that $s\geq 2$. The subset 
\begin{eqnarray*}
S&:=&
\{x\in \m Z/d\m Z\mid x\equiv 1 
\;{\rm mod}\; t\}\\
& =&\{ 1, 1\!+\!t, 1\!+\!2t,\ldots, 1\!+\!(s\!-\!1)t\}
\end{eqnarray*}
is a subgroup of $(\m Z/d\m Z)^*$ of order $s$,
and one has $xj\equiv j$ mod $d$, for all $x$ in $S$.
Since, by the second step, $k$ is coprime to $s$, the cardinality of the set $Sk$ is also equal to $s$ and one has
$$
Sk:=\{k'\in \m Z/d\m Z\mid k'\equiv k 
\;{\rm mod}\; t\}.
$$

When $r\geq 4$, since one has $s\geq 2$, the set $Sk$ contains an element $k'$ with $j\leq k'\leq d-j$.

When $r=3$, one has $s=3$ and the set $Sk$ still contains an element $k'$ with $j\leq k'\leq d-j$.
\vs 

{\bf Fifth step: When $j=2$, one has $\;3\,g(d)\geq d\!+\!2$.}\\
In this case, one has $k=1$. We remark that 
the interval $I:=[d/3,2d/3]$ does not contain integers coprime to $j$. 
Indeed, when an integer  $x$ is in this interval,
the integers $j':=|2x\!-\! d|$ and $k':=x$ satisfy 
\eqref{eqnjkdjkd}, i.e. one has the two inequalities
$$
j'\leq k'\leq d-j'.
$$
Therefore the integer $x$ is not coprime to $j$.
Since the number of integers in this interval $I$ is at least $(d-1)/3$, this proves the lower bound for the Jacobsthal function
\begin{equation}
\label{eqnlowbou1}
g(d)\geq (d\!+\!2)/3.
\end{equation}

{\bf Sixth step: One always has $\; r\, g(j) \geq (r\!-\!2)j\!+\! 2$.}\\
In this step, we will use again the fact that one can not find an integer $x$ coprime to $d$
for which the corresponding pair $(j',k')$ satisfies \eqref{eqnjkdjkd}.
We will only focus on the subgroup of $(\m Z/d\m Z)^*$
$$
S:=\{x\in (\m Z/d\m Z)^*\mid x\equiv 1\; {\rm mod}\; r\}.
$$ 
Note that, for all $x$ in $S$, the integer 
$j':=j$ satisfies $j'\equiv xj$ mod $d$.
Since $k$ and $r$ are both coprime to $j$, 
one can find integers $\ell$ and $s$ such that
$$
\ell k\equiv sr\equiv 1 \;{\rm mod}\; j.
$$

We will check that 
the interval 
$\displaystyle I:=\left[sk+\frac{j-k}{r},sk+\frac{d-j-k}{r}\right]$ 
does not contain integers coprime to $j$. 
Indeed, when an integealwaysr  $y=sk+m$ is in this interval and 
is coprime to  $j$,
the element 
$$
x:= 1+ m\ell r
\; \equiv\;
\ell r y
\;\;{\rm mod}\;\; j
$$
is also coprime to $j$, hence it belongs to $S$. Moreover the element 
$$
k':=k+(y\!-\!sk)r
\;=\;
k+mr 
\;\equiv\; 
xk\;{\rm mod}\; d
$$
satisfies
$$
j'\leq k'\leq d-j'.
$$
Therefore $y$ is not coprime to $j$.

Since the number of integers in this interval $I$ is at least
$(r-2)(j-1)/r$, this proves the lower bound for the Jacobsthal function
\begin{equation}
\label{eqnlowbou2}
 g(j) \geq( (r\!-\!2)j\!+\! 2)/r.
\end{equation}

{\bf Seventh step: We list the possible values of $j$ and $r$}.\\
We first notice that since $r\geq 3$, Inequality \eqref{eqnlowbou2} tells us that 
$g(j)>j/3$. Lemma \ref{lemkanold} tells us then that $j\leq 10$.
We then distinguish according to the values of $j$.
\vs

\noindent
{\bf $\star$ When $j=2$. } Inequality \eqref{eqnlowbou1} tells us that $g(d)> d/3$. Lemma \ref{lemkanold} tells us then that $d\leq 10$.
Since $r$ is coprime to $j$. This gives $r=3$ or $r=5$.
\vs 

\noindent
{\bf $\star$ When $j= 3,4,5,7,8$ or $9$. } In these cases one has $g(j)=2$. Inequality \eqref{eqnlowbou2} tells us that 
$(j\!-\!2)r\leq 2j\!-\!2$. Since $r\geq 3$, this implies 
that $j\leq 4$. 
Since $r$ is coprime to $j$, for $j=3$ this gives $r=4$, and for 
$j=4$ this gives  $r=3$.
\vs 

\noindent
{\bf $\star$ When $j= 6$ or $10$. } In these cases one has $g(j)=4$. Inequality \eqref{eqnlowbou2} tells us that 
$(j\!-\!4)r\leq 2j\!-\!2$. Since $r$ is coprime to $j$, for $j=6$ this gives $r=5$, and for 
$j=10$ this gives  $r=3$.
\vs

{\bf Eigth step: We analyze each of these six subcases}.
\vs 

\noindent
{\bf $\star$ When $j=2$ and $r=3$: }
One has $d=6$ and $k=1$. This is Case $(b)$.
\vs 

\noindent
{\bf $\star$ When $j=2$ and $r=5$: }
One has $d=10$ and $k=1$. This is Case $(c)$.
\vs 

\noindent
{\bf $\star$ When $j=3$ and $r=4$: }
One has $d=12$ and $k=1$ or $2$. 
The subcase $k=1$ is excluded since $(3,1)\sim (3,5)$.
The subcase $k=2$ is Case $(d)$.
\vs 

\noindent
{\bf $\star$ When $j=4$ and $r=3$: }
One has $d=12$ and $k=1$ or $3$. 
The subcase $k=1$ is excluded since $(4,1)\sim (4,7)$.
The subcase $k=3$ is Case $(e)$.
\vs 

\noindent
{\bf $\star$ When $j=6$ and $r=5$: }
One has $d=30$ and $k=1$ or $5$. 
The subcase $k=1$ is excluded since $(6,1)\sim (6,11)$.
The subcase $k=5$ is Case $(f)$.
\vs 

\noindent
{\bf $\star$ When $j=10$ and $r=3$: }
One has $d=30$ and $k=1, 3,7$ or $9$.\\ 
The cases $k\!=\!1$ or $7$ are excluded since $(10,1)\sim (10,7)\sim (10,19)$.\\
The cases $k\!=\!3$ or $9$ are equivalent since 
$(10,3)\!\sim\! (10,9)$. This is Case $(g)$.
\end{proof}

\subsection{The Stickelberger's formula}
\label{secstifor}

We will not use the full strength of Stickelberger's formula  that can be found in \cite[Chap.1]{LangCyclotomic} or
in \cite[Chap.14]{IrelandRosen}.
Instead we will use the following elementary formula
that was already known to Kummer and that is one of the key ingredients in its proof.

\bl
\label{lemstifor}
Let $\m F_p$ be a prime field, let $j,k$ be integers $0<j,k<p\!-\!1$. Let us define the Jacobi sum mod $p$ by 
$\; J_{j,k}:=\textstyle\sum\limits_{x\neq 0,1}x^{-j}(1-x)^{-k}
\in \m F_p$.\\
$a)$ One has the equality $J_{j,k}=
-\binom{j+k}{k}$ in $\m F_p$.\\
$b)$ In particular, one has the equivalence 
$J_{j,k}\neq 0 \Longleftrightarrow j+k<p$.
\el

Here the sum is over all $x$ in $\m F_p$ with $x\neq 0$, $x\neq 1$ and the right-hand side is the binomial coefficient
$\binom{j+k}{k}=\frac{(j+k)!}{j!\,k!} $.

\begin{proof}[Proof of Lemma \ref{lemstifor}]
 This is a classical and  elementary calculation
 \begin{eqnarray*}
 J_{j,k}&=&\textstyle\sum\limits_{x\neq 0}x^{-j}(1-x)^{p-1-k}
 \;=\;
 \textstyle\sum_{\ell=0}^{p-1-k}(-1)^\ell\binom{p-1-k}{\ell}\sum\limits_{x\neq 0}
 x^{\ell-j}\\
 &=&\textstyle -(-1)^j\binom{p-1-k}{j}
 \;=\;
 -\binom{j+k}{k},
 \end{eqnarray*}
which is valid since
the base field is $\m F_p$.
\end{proof}

\subsection{Distinguishing Jacobi sums}
\label{secdisjac}

We can now give the proof of $\Longrightarrow$ in Proposition \ref{projacsum}.
\vs 

\noindent
{\bf Notation } 
We will prove Proposition \ref{projacsum} by reducing 
\eqref{eqnratjac} modulo a suitable prime ideal  $\g p$ of the ring $\m Z[\zeta_{p-1}]$. More precisely, the multiplicative group $\m F_p^*$ is a cyclic group
of order $d:=p\!-\! 1$.
We denote by $g_0$ the smallest positive integer whose image modulo $p$ is a generator
of $\m F_p^*$
and we denote by $\om$ the Teichm\"{u}ller character of $\m F_p^*$
which is defined by the equality $\om(g_0)=\zeta_{p-1}$.
This character is a generator of the group of characters of $\m F_p^*$. 
In particular, since $\chi_1$ and $\chi_2$ are not principal, there exist positive integers 
$j,k<p\!-\! 1$ such that
$$\chi_1=\om^{-j}
\; ,\;\;
\chi_2=\om^{-k}
\; ,\; \mbox{\rm and hence }\;\;
\ol{\chi}_1=\om^{-(p-1-j)}.
$$

\begin{proof}[Proof of the implication $\Longrightarrow$ in Proposition \ref{projacsum}]

We assume that we are not in one of the seven cases $(a)\!-\!(g)$. In particular,  one has $j\neq \frac{p-1}{2}$.
The action of an element of the Galois group of $K/\m Q$
commutes with the complex conjugation and hence 
preserves the assertion \eqref{eqnratjac}.
This action 
is given by an element $a\in (\m Z/d\m Z)^*$ and 
sends the characters $\om^{-j}$ and $\om^{-k}$ respectively to 
the characters $\om^{-aj}$ and $\om^{-ak}$. 
Therefore, without loss of generality by using the combinatorial Lemma \ref{lemcomzdz}, 
we can assume that
$$
j\leq k\;\;{\rm  and }\;\; j+k\leq p\!-\! 1.
$$

We first deal with the case where  $j=k$. In this case, one has
$j=k<\frac{p-1}{2}$ and, by \eqref{eqnj12srp}, the Jacobi sums in \eqref{eqnratjac} do not have the same absolute value.
This proves that the ratio \eqref{eqnratjac}
is not a root of unity.

We now deal with the case where $j<k$. The Jacobi sum $J(\chi_1,\chi_2)$ lives in the ring of integers ${\mc O}_K:=\m Z[\zeta_{p-1}]$. 
Since the polynomial $X^{p-1}-1$ has $p\!-\! 1$ distinct root 
in $\m F_p$, the cyclotomic polynomial $\Ph_{p-1}(X)$ is also 
split in $\m F_p$ and has $\ph(p\!-\!1)$ roots in $\m F_p$ where $\ph$ is the Euler totient. These roots are the generators of the group $\m F_p^*$. We denote by $\g p$ the prime ideal of $\mc O_K$ 
over $p$ containing $g_0\!-\!\zeta_{p-1}$. We denote by 
$$
\pi_{\g p}: \mc O_K\ra
\mc O_K/\g p\simeq \m F_p
$$
the natural morphism given by the reduction modulo $\g p$ so that one has 
$\pi_{\g p}(\om(x))=x$ for all $x$ in $\m F_p$. 
Using the notation and the result of Lemma \ref{lemstifor}, one has the equalities
$$
\pi_{\g p}(J(\chi_1,\chi_2))=\pi_{\g p}(J(\om^{-j},\om^{-k}))=J_{j,k}=
\textstyle \binom{j+k}{k}\neq 0
$$
since $j\!+\!k \leq p\!-\! 1$. Similarly, one has the equalities
$$
\pi_{\g p}(J(\ol{\chi}_1,\chi_2))=\pi_{\g p}(J(\om^{-(p-1-j)},\om^{-k}))=J_{p-1-j,k}=
\textstyle \binom{p-1-j+k}{k}= 0
$$
since $j<k$. This proves that the ratio \eqref{eqnratjac}
is not a root of unity.
\end{proof}

\subsection{Example of equalities between Jacobi sums}
\label{secexajac}

In this section we give the proof of $\Longleftarrow$ in
Proposition \ref{projacsum}.
We also explain that the ratios \eqref{eqnratjac}
are frequently equal to $1$ in each of the seven cases.

\begin{proof}[Proof of $\Longleftarrow$ in 
Proposition  \ref{projacsum}] 
We want to prove that in each of the seven cases of Proposition \ref{projacsum},
the ratio  $R(\chi_1,\chi_2)$ in \eqref{eqnratjac}
is a $(p-1)^{\rm th}$-root of unity.

This ratio is the ratio of two algebraic integers
in the cyclotomic field $K=\m Q[\zeta_{p-1}]$ all of whose conjugate have the same  absolute value equal to $\sqrt{p}$.
The pairs $(\chi_1,\chi_2)=(\omega^{-j},\omega^{-k})$, 
in these seven cases 
are exactly those for which one can not find a pair $(j',k')$ satisfying \eqref{eqnjkdjkd}.
By Stickelberger's formula,this means that the set of prime 
ideals $\g p$ 
of $\mc O_K$
over $p$ that contain $J(\chi_1,\chi_2)$ is exactly the same as the set of prime ideals of $\mc O_K$
over $p$ that contain $J(\ol{\chi}_1,\chi_2)$. 
Since $J(\chi_1,\chi_2)$ has absolute value $\sqrt{p}$, 
and since the ideal $(p)$ of $\mc O_K$ is a product of distinct prime ideals,  the ratio $R(\chi_1,\chi_2)$ is a unit in $K$. 
Since all the Galois conjugates of this ratio $R(\chi_1,\chi_2)$ have absolute value equal to $1$, this ratio  is a root of unity in $K=\m Q[\zeta_{p-1}]$, hence it is a
$(p-1)^{\rm th}$-root of unity.
\end{proof}

\br
\label{remcasbcg} One can also check the following facts that we will not use.\\
In Cases $(a)$ and $(e)$, the ratio 
$R(\chi_1,\chi_2)$ is always equal to $1$.\\
In Case $(b)$ and $(g)$, the ratio $R(\chi_1,\chi_2)$ is always a cubic root of unity.\\
In Cases $(c)$ and  $(f)$, the ratio $R(\chi_1,\chi_2)$ is always a fifth root of unity.\\
In Case $(d)$, the ratio $R(\chi_1,\chi_2)$ is always equal to $\pm 1$.
\er

\noindent{\bf Examples}:
In these seven cases, one does not always have the equality
$R(\chi_1,\chi_2)=\chi_2(-1)$. This depends on the prime $p$ and on the parity of $\chi_2$. But in each of these seven cases, there are many primes $p$ for which 
\begin{equation}
\label{eqnrc1c21}	
R(\chi_1,\chi_2)=\chi_2(-1)=1.
\end{equation}
One can compute the smallest value of $p$ for which 
\eqref{eqnrc1c21}	 holds:\\
In Case $(a)$: $p= 7$. One has 
$J(\om^{3},\om^{2})= J(\om^{-3},\om^{2}) =2-i\sqrt{3}$.\\
In Case $(b)$: $p\!=\! 109$. One has 
$J(\om^{36},\om^{18})\!=\! J(\om^{-36},\om^{18})\! =-1-6i\sqrt{3}$.\\
In Case $(c)$: $p=241$. One has 
$J(\om^{48},\om^{24})= J(\om^{-48},\om^{24})$.\\
In Case $(d)$: $p=37$. One has 
$J(\om^{9},\om^{6})= J(\om^{-9},\om^{6})=-5+2i\sqrt{3}$.\\
In Case $(e)$: $p=73$. One has 
$J(\om^{24},\om^{18})= J(\om^{-24},\om^{18})=8+3i$.\\
In Case $(f)$: $p=601$. One has 
$J(\om^{120},\om^{100})= J(\om^{-120},\om^{100})=\frac{23-25i\sqrt{3}}{2}$.\\
In Case $(g)$: $p=421$. One has 
$J(\om^{140},\om^{42})= J(\om^{-140},\om^{42})$. 

\br
\label{remcasbcg0} In these seven cases, when $R(\chi_1,\chi_2)=1$ and $\chi_2$ is even,
one can also check the following facts that we will not use.\\
In Case $(a)$, $J(\chi_1,\chi_2)$ 
belongs to $\m Z[e^{2 i\pi/(p-1)}]$.\\
In Cases $(b)$, $(d)$ and $(f)$, $J(\chi_1,\chi_2)$ 
belongs to $\m Z[e^{2 i\pi/3}]$.\\
In Case $(c)$ and $(g)$,  $J(\chi_1,\chi_2)$ 
belongs to $\m Z[e^{2 i\pi/5}]$.\\
In Case $(e)$, $J(\chi_1,\chi_2)$ 
belongs to $\m Z[e^{2 i\pi/4}]$.
\er

This ends the proof of Proposition \ref{projacsum}   and hence this also ends 
the proof of Theorem \ref{thmexicoh} for non safe primes.

\section{Unimodular $\mc C$-functions for safe primes}
\label{secequsaf}

The aim of this chapter is to prove Theorem \ref{thmexicoh}
when the prime $p$ is safe, i.e. when the integer $n:=(p-1)/2$ is prime.
We will work with the projective space $\m P(V)$ of the space $V$ of odd functions on $\m F_p$, its Clifford torus $T$ and the Fourier transform $F$.
The new difficulty in this case is that the  Clifford tori $T$ and $F^{-1}T$ are not transverse at the intersection point $[\chi_0]$
given by the Legendre character $\chi_0$. Indeed they are tangent at that point.

To overcome this difficulty, we will study the ``Hessian'' of this intersection and prove that it is anisotropic in Lemmas \ref{lemhesani} and \ref{lemasupsi}.
This will prove that the algebraic multiplicity at $[\chi_0]$ of the intersection $T\cap F^{-1}T$ is equal to $2^{n-1}$, see 
Proposition \ref{protrach0}.$a$.

At first glance, since this number $2^{n-1}$ is the one predicted by 
Theorem \ref{thmchobis}, this information seems useless.
To go on we will find a small unitary perturbation $F_t$ of the Fourier transform $F$ so that, in a neighborhood of $[\chi_0]$ 
the intersection $T\cap F_t^{-1}T$ is transverse and contains at most $2^{n-2}$ points, see Proposition \ref{protrach0}.$b$. This will allow us to conclude.

\subsection{Strategy}
\label{secstrsaf}

In this section we explain in more detail the strategy 
and prove Theorem \ref{thmexicoh} for safe primes,
taking for granted Proposition  \ref{protrach0} that will be proven later in this chapter.
We will keep the following notation 
for most of Chapter \ref{secequsaf}.\vs 

Let $p\geq 7$ be a safe prime. Let $H_0=\{\pm 1\}\subset \m F_p^*$ and let $c_0$ be the non trivial character of $H_0$. Let 
\begin{equation*}
V=V_{H_0,c_0}=\{ f: \m F_p \ra \m C\mid \; f \;\mbox{\rm is odd}\,\},
\end{equation*}
\begin{equation*}
T=T_{H_0,c_0}=\{ [f]\in \m P(V)\mid \; |f| \;\mbox{\rm is constant}\,\}. 
\end{equation*}
Let $F:V\ra V$ be the Fourier transform, let $\chi_0\in V$ be the Legendre character and $n=(p-1)/2=\dim V$. By assumption, the integer $n$ is an odd prime number. The space $V$ is a hermitian vector space for the $\ell^2$ norm and $F$ is a unitary transformation of $V$. 

The following proposition is analogous to Proposition \ref{procouni}.$b$.

\bp
\label{protrachi} Let $p$ be a prime with $p\equiv 3$ mod $4$.
For all odd Dirichlet character $\chi\neq \chi_0$, the intersection
$T\cap F^{-1}T$ is transverse at $[\chi]$.
\ep

\begin{proof}[Proof of Proposition \ref{protrachi}]
Assume that this intersection $T\cap F^{-1}T$ is not transverse at $[\chi]$. We apply Propositions \ref{projacsum} and \ref{protracli}.
We note that, since the index $n=(p\!-\!1)/2$
of $H_0$ is odd, all the characters $\chi_2=\psi$ of $\m F_p^*/H_0$ have odd order. Therefore the only option among the seven cases of Proposition \ref{projacsum} is the first one: the character $\chi$ must have order $2$. Hence this character is the Legendre character $\chi_0$.
\end{proof}

We recall that 
$T$ and $T':=F^{-1}T$ are smooth real algebraic subvarieties of $\m P(V)$ and that the 
algebraic multiplicity at $[\chi_0]$ of the intersection $T\cap T'$ is by definition the dimension 
\begin{equation}
\label{eqnmulint}
 \mbox{$m_0:=\dim _{\m C}A_0$. Here $A_0$ is  the local ring 
$A_{0}:=\mc A_0/(\mc I_0\cap \mc I'_0)$}
\end{equation}
 where $\mc A_0$ is the local ring
at $[\chi_0]$ of germs of complex valued real analytic functions on  $\m P(V)$
seen as a real analytic variety, and  where $\mc I_0$ and $\mc I'_0$ are respectively
the ideals of germs of functions that vanish on $T$ and on $T'$ respectively. 

\bp
\label{protrach0}
Let $p\geq 7$ be a safe prime.\\ 
$a)$ The tori $T$ and $F^{-1} T$ are tangent at $[\chi_0]$
and the algebraic multiplicity at $[\chi_0]$ of the intersection 
$T\cap F^{-1}T$ is $2^{n-1}$.\\
$b)$ There exist a neighborhood $\Om_{\chi_{_0}}$ of $[\chi_0]$ in $\m P(V)$ and  a continuous family of unitary maps $F_t$ of $V$ with $F_0=F$ such that, for all $t\in (0,1)$ the intersection $T\cap F_t^{-1}T$ is transverse in $\Om_{\chi_{_0}}$ and contains at most $ 2^{n-2}$ points in $\Om_{\chi_{_0}}$.
\ep

Note that the upper bound of the number of intersection points near $[\chi_0]$ for the perturbation of the Fourier transform is  half the algebraic  multiplicity at $[\chi_0]$ of this intersection for the Fourier transform.

\begin{proof}[Proof of Proposition \ref{protrach0}] It will occupy the rest
of Chapter \ref{secequsaf}. 

Proposition \ref{protrach0}.$a$ 
will follow from  Lemma \ref{lemasupsi}. 
 
Proposition \ref{protrach0}.$b$
will follow from 
Lemma  \ref{lemmulfam}.
\end{proof} 

\begin{proof}[Proof of Theorem \ref{thmexicoh} for a safe prime]
Assume, by contradiction that the intersection $T\cap F^{-1}T$ contains only
the $n$ points $[\chi]$ associated to the $n$ Dirichlet characters $\chi$ that belong to $V$.
We introduce an open neighborhood $\Om_{\chi_{_0}}$ of $[\chi_0]$ and a unitary perturbation $F_t$ of the Fourier transform $F$ as in Proposition \ref{protrach0}. We know that, for all $t\in (0,1)$, 
the intersection $T\cap F_t^{-1}T$ is transverse in $\Om_{\chi_{_0}}$ and contains at most $ 2^{n-2}$ points in $\Om_{\chi_{_0}}$. 

By Proposition \ref{protrachi}  the intersection $T\cap F^{-1}T$ is transverse at all 
the points $[\chi]$ when $\chi\neq \chi_0$. 
Hence, for $\chi\neq \chi_0$, there also exists a neighborhood 
$\Om_\chi$ of $[\chi]$ in $\m P(V)$ such that, if $t$ is small enough,
the intersection  $T\cap F_t^{-1}T$ is transverse in $\Om_\chi$ and contains at most one point in $\Om_\chi$.

One moreover choose $t$ small enough so that the intersection  $T\cap F_t^{-1}T$ is included in the union of all the neighborhoods $\Om_\chi$.
Then this intersection is transverse and contains at mos $2^{n-2}+n-1$ points. Theorem \ref{thmchobis} predicts the existence of 
at least $2^{n-1}$ points in this intersection.
Since $n\geq 5$, one has $2^{n-1}>2^{n-2}+n-1$. 
This is the contradiction we are looking for.
\end{proof}

\subsection{A quadratic functional equation}
\label{secquaequ}

We begin by a crucial lemma that solves a quadratic equation
and will allow us to compute the multiplicity in Proposition \ref{protrach0}.$a$. 

\bl
\label{lemhesani} 
Let $p$ be a safe prime and $\al: \m F_p\ra \m C$ be an odd function with\\
$(i)$ \; $\al *\al =\chi_0 * \chi_0\al^2$,\\
$(ii)$ \; $\langle \chi_0,\al\rangle=0$.\\
Then one has $\al=0$.
\el
Note that if $\al$ is an odd function solution of $(i)$, then 
for all $t$ in $\m C$, the odd function $\al +t\chi_0$ is also an  odd function solution of $(i)$. This is why the linear condition $(ii)$ 
is natural and useful.

Note also that Equation $(i)$ can be rewritten more explicitely as 
\begin{equation}
\label{eqnhesanj}
\textstyle 
\sum\limits_{\ell\in \m F_p}\al_{k-\ell}\al_\ell
\; =\; 
\sum\limits_{\ell\in \m F_p} \,
(\!\frac{k-\ell}{p}\!)(\!\frac{\ell}{p}\!)\,\al_\ell^2
\;\;\;\;\mbox{\rm for all $k$ in $\m F_p$}.
\end{equation}

Note finally that Equation $(ii)$ can be rewritten as, 
\begin{equation}
\label{eqnhesank}
\textstyle 
\sum\limits_{\ell\in C}\,\al_\ell
\; =\; 0.
\end{equation}
where $C=(\m F_p^*)^2$ is the 
set of non zero squares. Indeed, since $p\equiv 3$ mod $4$, one has $\m F_p^*=C\cup -C$.

\begin{proof}
The main trick in the proof is a reduction modulo $2$ argument.
We first identify the field $\m C$ of complex numbers with the algebraic closure $K:=\ol{\m Q}_2$ of the local field $\m Q_2$. We denote by 
$v_2:K\ra \m Q\cup\{+\infty\}$ the $2$-adic valuation and by $\mc O_K=\{\la \in K\mid v_2(\la)\geq 0\}$ the ring of integers of $K$. It is a local ring  whose maximal ideal  is $\g m_K:=\{\la \in K\mid v_2(\la)> 0\}$ and whose residual field $\ka=\mc O_K/\g m_K$ is the algebraic closure 
$\ol{\m F}_2$ of $\m F_2$.

Let $\al:\m F_p\ra\ol{\m Q}_2$ be an odd function 
satisfying $(i)$ and $(ii)$. We want to prove that $\al=0$.
Assume by contradiction that $\al$ is non zero. 
Replacing $\al$ by a multiple, one can assume that $\al$ takes its values in $\mc O_K$ and that the reduction 
$\dt{\al}:\m F_p\ra \ol{\m F}_2$ of $\al$ modulo $\g m_K$ is non-zero.
Since $\al$ is odd, one  has 
\begin{equation}
\label{eqnhesmoi}
\dt{\al}_0\; =\; 0 \;\;\;{\rm and}\;\;\;
\dt{\al}_k\; =\; \dt{\al}_{-k}
\;\;\;\;\mbox{\rm for all $k$ in $\m F_p$}.
\end{equation}
We write Equation \eqref{eqnhesanj}  with $2k$ instead of $k$ and reduce it modulo $2$.
We notice that all  the terms in the sums occur by pairs, except for $\ell=k$ in the left-hand side and except for $\ell=-2k$ in the right-hand side. This gives
\begin{equation*}
\dt{\al}_k^2\; =\; \dt{\al}_{2k}^2
\;\;\;\;\mbox{\rm for all $k$ in $\m F_p$}.
\end{equation*}
Since the function $\dt\al$ takes its values in a field of characteristic $2$ this implies
\begin{equation}
\label{eqnhesmoj}
\dt{\al}_k\; =\; \dt{\al}_{2k}
\;\;\;\;\mbox{\rm for all $k$ in $\m F_p$}.
\end{equation}
Reducing Equation \eqref{eqnhesank} modulo $2$ gives
\begin{equation}
\label{eqnhesmok}
\textstyle\sum\limits_{\ell\in C}\dt{\al}_\ell\; =\; 0.
\end{equation}
Since the integer $(p-1)/2$ is prime, the group $\m F_p^*$ is generated by the two elements 
$-1$ and $2$. Therefore, Equations \eqref{eqnhesmoi} and \eqref{eqnhesmoj} tell us that the function $\dt\al$ is constant equal to $\dt\al_1$ outside $0$. Since the cardinality of $C$ is odd, plugging this information into \eqref{eqnhesmok} gives $\dt\al_1=0$, and therefore $\dt\al=0$. This is a contradiction.
\end{proof}


\subsection{Multiplicity of holomorphic maps}
\label{secmulhol}


In this section, we recall a few general definitions and facts that deal with the topological degree of a proper holomorphic map with finite fibers. 
\vs 

For $\Om\subset \m C^m$ an open set, $\Psi_0:\Om\ra \m C^m$
a holomorphic map  and $z\in \Om$ such that $\Psi_0(z)=0$, 
we recall that  the {\it algebraic multiplicity} $m_z(\Psi_0)$
is the dimension of the local ring $\mc O_z/\Psi_0^*\mc O_{0}$
where $\mc O_z$ is the local ring of germs at $z$ of holomorphic functions
and where $\Psi_0^*\mc O_{0}$ is the ideal spanned by the functions $h\circ \Psi_0$
where $h$ is a holomorphic function that vanishes at $0$.
A point $z_0\in \Om$ is a {\it critical point} for $\Psi_0$ if and only if the tangent map
$D\Psi_0(z_0)$ is not invertible, or equivalently $m_{z_0}(\Psi_0-f(z_0))>1$. A {\it critical value} for $\Psi_0$ is the image 
$\Psi_0(z_0)$ of a critical point $z_0$. A {\it regular value} is a value which is not critical.
By Sard theorem, the set of critical values has Lebesgue measure $0$.
When the algebraic multiplicity $m_z(\Psi_0)$ is finite it coincides with 
the {\it geometric multiplicity}. This means that, there exists $\eps_0>0$ and a 
neighborhood $\Om_{0}$ of $z$ such that for Lebesgue almost all 
$w\in \m C^m$ with $\| w\|<\eps_0$, the value $w$ is a regular value of $\Psi_0\!-\! w$ and one has $m_z(\Psi_0)=\#(\Psi_0^{-1}(w)\cap \Om_{0})$.  See \cite[Ch.1 prop. 2.1]{AizenbergYuzhakov} and \cite[p.148]{Tsikh}.
\vs 

The first fact extends to several variables the famous Rouch\'e theorem.

\bfa
\label{facmulhol} 
Let $\Om\subset \m C^m$ be an open set, let ${\ol B}\subset \Om$ be a compact ball and $I\subset \m R$ be an interval. Let $(\Psi_t)_{t\in I}$ be a continuous 
family of holomorphic map $\Psi_t:\Om\ra \m C^m$ such that for all $t$ in $I$, 
the function $\Psi_t$ does not vanish on the boundary $\partial B$.
Then the number of zeros of $\Psi_t$ in $\ol B$ counted with multiplicities :
\;\;
$\sum\limits_{\{ z\in \ol B\mid\Psi_t(z)=0\} } m_z(\Psi_t)$\;\; 
does not depend on $t$.
\efa

\begin{proof}
See \cite[Thm 2.5]{AizenbergYuzhakov}.
\end{proof}

The second fact explains how to compute the algebraic multiplicity
for a complete intersection quadratic map. It also explains why in some fibers at most half of the points are real.

\bfa
\label{facintqua} 
Let $Q=(Q_1,\ldots,Q_m):\m C^m\ra \m C^m$ be a homogeneous polynomial map of degree $2$ which is anisotropic over $\m C$ i.e. such that $Q^{-1}(0)=\{0\}$.\\ 
$a)$ Then there exists $C>0$ such that, for all $z\in \m C^m$, one has
\begin{equation*}
\label{eqnintqua}
C^{-1}\| z\|^2\leq \|Q(z)\| \leq C\| z\|^2.
\end{equation*}
$b)$ The algebraic multiplicity of $Q$ at $0$ is $m_0(Q)=2^m$.\\
$c)$ For all $w\neq 0$ in $\m C^m$, there exists a sign $\eps=\pm 1$ 
such that 
\begin{equation*}
\label{eqnfibrea}
\# (Q^{-1}(\eps w)\cap \m R^m)\; \leq \; 2^{m-1}. 
\end{equation*}
\efa

\begin{proof}
$a)$ These inequalities follow from the homogeneity of $Q$ by choosing the constant $C$ to be  $C=\max\limits_{\|z\|=1}\max(\| Q(z)\|,\| Q(z)\|^{-1})$.

$b)$ According to \cite[p.142]{Tsikh}, since $Q$ is anisotropic, the multiplicity of $Q$ at $0$ is equal to the product of the degrees of the homogeneous polynomials $Q_j$ which is $2^m$.

$c)$ Remember the equality $Q(iz)=-Q(z)$, for all $z$ in $\m C^m$. For $w\neq 0$, the set $Q^{-1}(w)\cup Q^{-1}(-w)$ contains at most $2^{m+1}$ points and is stable by multiplication by $i$. Therefore it contains at most $2^m$ real points and one of these two fibers $Q^{-1}(\eps\, w)$ contains at most $2^{m-1}$ real points. 
\end{proof}

\br We will not need it, but one can say much more: the 
quotient algebra $A_0:=\m C[z_1, ,z_m]/\langle Q_1,\ldots , Q_m\rangle$ is a finite graded commutative
$\m C$-algebra
whose Hilbert polynomial is $P_{A_{_0}}(t)=(1+t)^m$. 
Moreover this local algebra $A_0$ is a ``complete intersection algebra'', hence a
``Poincaré duality algebra'' and a ``Cohen-Macauley algebra''.
These facts follow from the Macaulay theorem and the Koszul theorem (see \cite[thm 6.7.6 and 6.2.4]{SmithInvariant} ).
\er

The third fact tells us that one can sometimes reduce the 
calculation of the algebraic multiplicity to the previous example.

\bfa 
\label{facmulpsi}
Let $\Om\subset \m C^m$ be an open set containing $0$ and $\Psi_0:\Om\ra \m C^m$ be a holomorphic map such that  $\Psi_0(0)=0$ and $D\Psi_0(0)=0$.
Assume that the Hessian $Q:=\frac12 D^2\Psi_0(0):\m C^m\ra \m C^m$ is an anisotropic quadratic map.\\
$a)$ Then there exist $\eps_0>0$,  $C>0$ such that, for $z\in \m C^m$ with $\|z\|\leq \eps_0$, 
\begin{equation*}
	\label{eqnintpsi}
	C^{-1}\| z\|^2\leq \|\Psi_0(z)\| \leq C\| z\|^2.
\end{equation*}
$b)$ The algebraic multiplicity of $\Psi_0$ at $0$ is $m_0(\Psi_0)=2^m$.
\efa

\begin{proof}
$a)$ We write $\Psi_0(z)=Q(z)+R(z)$ where $Q:\m C^m\ra \m C^m$ is a 
quadratic polynomial map and where $\lim\limits_{z\ra 0}\frac{\|R(z)\|}{\|z\|^2}=0$. We apply then Fact \eqref{facintqua}.$a$.
	
$b)$ According to \cite[p.146]{Tsikh}, since $Q$ is anisotropic, the multiplicity of $\Psi_0$ at $0$ is equal to the multiplicity of $Q$ at $0$
which is $2^{m}$.
\end{proof}

\subsection{Multiplicity of $T\cap F^{-1}T$ at the Legendre character}
\label{secmulleg}

We come back to the notation of Sections \ref{secstrsaf} and \ref{secquaequ}. 
The aim of this section is to apply the general results of Section \ref{secmulhol} to determine the multiplicity of the intersection 
$T\cap F^{-1}T$ at the point $[\chi_0]$.
\vs 

We introduce the complex vector spaces
\begin{eqnarray*}
V_0
&=& 
\{\al: \m F_p\ra \m C\mid
\mbox{\rm $\al$ is odd and 
$\langle \chi_0,\al\rangle =0$} \},\\
W_0
&=& 
\{\be: \m F_p\ra \m C\mid
\mbox{\rm $\be$ is even and 
	$\be(0)=\widehat{\be}(0) =0$} \}.
\end{eqnarray*}

Note that these two vector spaces have dimension $n\!-\! 1=(p\!-\! 3)/2$ and that
the formula $\al=\chi_0\be$ induces an isomorphism between $V_0$ and $W_0$.

It will be convenient to use a parametrization of a neighborhood of $[\chi_0]$ in $\m P(V)$ by the ball 
\begin{eqnarray*}
B_0
&=& 
\{\be\in W_0\mid\;
\| \be\|_\infty <\pi/2 \}.
\end{eqnarray*}
This parametrization is given by 
\begin{equation}
\label{eqnbetfbe}
{\be}\mapsto [f_{\be}]
\;\;{\rm where}\;\; 
f_{\be}=\chi_0 \,e^{i\be}\, .
\end{equation}
The next lemma describes $T$ and $F^{-1}T$ 
via this parametrization.
We set
\begin{equation*}
\label{eqnvowobo}
V_{0,\m R}\!=\!\{\al\!\in\! V_0\mid \al\!=\!\ol{\al}\}\, ,\,
W_{0,\m R}\!=\!\{\be\!\in\! W_0\mid \be\!=\!\ol{\be}\}
\;{\rm and}\;
B_{0,\m R}\!=\! B_0\!\cap\! W_{0,\m R}.
\end{equation*}

\bl
\label{lemmulmul} 
$a)$ The following  map $\Psi_0: W_0\ra W_0$ is well defined
\begin{equation}
\label{eqndefpsi}
{\be}\mapsto \Psi_0(\be):=\,
\chi_0e^{i\be}\, *\, \chi_0e^{-i\be}\; -\; \chi_0 *\chi_0\, .
\end{equation}
$b)$ For $\be$ in $B_0$, one has the equivalences:
$\; \chi_0e^{i\be}\in T
\Longleftrightarrow \be=\ol\be\; \;$ and\\
$\; \chi_0e^{i\be}\in F^{-1}T
\Longleftrightarrow 
\chi_0e^{i\be}\, *\, \chi_0e^{-i\ol\be}\; =\; \chi_0 *\chi_0.$\\
$c)$ The multiplicity $m_0$ of the intersection $T\cap F^{-1}T$ at 
$[\chi_0]$ is equal to the algebraic multiplicity $m_0(\Psi_0)$ of $\Psi_0$ at $0$.
\el

Note that, since $p\equiv 3$ mod $4$, one has  $\chi_0 *\chi_0= {\bf 1}_{\m F_p}-p\,\de_0$.

\begin{proof}
$a)$ We want to check that, for $\be$ in $W_0$, the function $\Psi_0(\be)$ also belongs to $W_0$.
Since $p\equiv 3$ mod $4$, the Legendre character $\chi_0$ is odd and the 
functions $\chi_0 e^{\pm i\be}$ are also odd. 
Therefore the function $\Psi_0(\be)$ is even.\\
Since the function $\be$ is even, one has
$$
\Psi_0(\be)(0)= \sum_{k\in \m F_p}\;(\chi_0(k)e^{i\be(k)}\chi_0(-k)e^{-i\be(-k)}- \chi_0(k)\chi_0(-k))=0
$$
Since  the functions $f_{\be}=\chi_0 e^{ i\be}$ is odd, one has
$$\textstyle
\frac{1}{\sqrt{p}}\,\widehat{\Psi_0(\be)}(0)=\widehat{f_\be}(0)\;\widehat{f_{-\be}}(0)- 
(\widehat{\chi_0}(0))^2=0
$$
This proves that the function $\Psi_0(\be)$ belongs to $W_0$.

$b)$ The first equivalence follows from the definition. The proof of the second equivalence 
is similar to the proof of Lemma \ref{lemfoucohn}. 

$c)$ In the coordinate system $\be\mapsto [f_\be]$
as in \eqref{eqnbetfbe}, the components of the real analytic map $\be\mapsto \be-\ol{\be}$ 
spans the ideal of germs at $[\chi_0]$ of real analytic functions that vanish on $T$.
Similarly, the components of the real analytic map $\be\mapsto \chi_0e^{i\be}\, *\, \chi_0e^{-i\ol\be}- \chi_0 *\chi_0$ 
spans the ideal of germs at $0$ of real analytic functions that vanish on $F^{-1}T$. 
Therefore the local ring $A_0$ at $[\chi_0]$ of the schematic intersection $T\cap F^{-1}T$ defined in \eqref{eqnmulint} is isomorphic to $\mc O_0/\Psi_0^*(\mc O_0)$
where $\mc O_0$ is the ring of germs at $0$ of holomorphic functions on $W_0$. The dimension of this ring is $m_0(\Psi_0)$.
\end{proof}

We check now that this map $\Psi_0$ satisfies the assumptions
of Fact \ref{facmulpsi}

\bl 
\label{lemasupsi}
The holomorphic  map $\Psi_0$ given by \eqref{eqndefpsi} satisfies:\\
$(i)$ at the point $0$, the map $\Psi_0$ vanishes: $\Psi_0(0)=0$.\\
$(ii)$ the differential $D\Psi_0(0): W_0\ra W_0$ vanishes: $D\Psi_0(0)=0$,\\
$(iii)$ the hessian $Q=\frac12 D^2\Psi_0(0): W_0\ra W_0$ is an anisotropic  quadratic map.\\
And therefore the muliplicity $m_0=m_0(\Psi_0)$ is equal to $2^{n-1}$.
\el 

\br 
The assertion $(ii)$ reflects the fact that the Clifford tori $T$ and $F^{-1}T$
are tangent at the point $[\chi_0]$.
\er

\begin{proof}
One computes, for $\be$ in $W_0$, using Formula \eqref{eqndefpsi}:
\begin{eqnarray*}
\Psi_0(0)&=&\chi_0 *\chi_0-\chi_0*\chi_0=0,\\
D\Psi_0(0)(\be)&=&i\chi_0 \be *\chi_0 -i\chi_0 * \chi_0\be=0,\\
Q(\be)=\tfrac12 D^2\Psi_0(0)(\be)&=&\chi_0\be *\chi_0\be -\chi_0* \chi_0\be^2.
\end{eqnarray*}
According to Lemma \ref{lemhesani} this quadratic map $Q$ is anisotropic.
Therefore, we can apply Fact \ref{facmulpsi} to the holomorphic map $\Psi_0$
and the multiplicity of $\Psi_0$
at $0$ is equal to $2^m$ where $m=\dim_\m C W_0=n\!-\!1=(p\!-\!3)/2$.
\end{proof}

\subsection{Moving slightly the torus}
\label{secmovtor}

The aim of this section is to prove Proposition \ref{protrach0}, that is to  find a small unitary deformation
$F_t$ of the Fourier transform $F$ such that the intersection 
$T\cap F_t^{-1}T$ is transverse and contains at most $2^{n-2}$ points.

In order to clarify our arguments, we begin by
a general result that deal with analytic families of real analytic maps.

\bl
\label{lemmuldef} 
Let $I\subset \m R$ and $\Om\subset \m C^m$ be two open sets  containing $0$, and let $\Psi:I\times \Om\ra \m C^m; (t,z)\mapsto \Psi_t(z)$ be a real analytic family of holomorphic maps such that\\
$(i)$ the value $\Psi_0(0)$ vanishes: $\Psi_0(0)=0$,\\
$(ii)$ the differential  $D\Psi_0(0): \m C^m\ra \m C^m$ vanishes: 
$D\Psi_0(0)=0$,\\
$(iii)$ the hessian $Q:=\tfrac{1}{2}D^2\Psi_0(0): \m C^m\ra \m C^m$ is anisotropic: $Q^{-1}(0)=\{0\}$,\\
$(iv)$ the element $w_0:=-\tfrac{d}{dt}\Psi_t(0)|_{t=0}\in \m C^m$ 
is a regular value of $Q$.\\
Then there exists $t_0>0$, a radius $\rho_0>0$ and a sign $\eps=\pm 1$  such that, for all $t$ in $(0,t_0)$, 
the point $0$ is a regular value for the restriction of $\Psi_{\eps t}$ 
to the ball $B_{\m R^m}(0,\rho_0)$ and one has 
\begin{equation}
\label{eqnpstpsm}
\#(\Psi_{\eps t}^{-1}(0)\cap  B_{\m R^m}(0,\rho_{0}))
\;\;\leq\;\;  2^{m-1}.
\end{equation}
\el

Here $B_{\m R^m}(0,\rho_0)$ is the ball in $\m R^m$ with center $0$ and radius $\rho_0$, and the bound \eqref{eqnpstpsm} is a control on the number of real points near $0$
in the fiber $\Psi_{\eps t}^{-1}(0)$ by half the multiplicity of $\Psi_0$ at $0$.

\begin{proof}
According to Fact \ref{facintqua}, there exists a sign $\eps=\pm 1$ 
such that
\begin{equation}
\label{eqnfibrea2}
\# (Q^{-1}(\eps \, w_0)\cap \m R^m)\; \leq \; 2^{m-1}. 
\end{equation}
Note that, since $Q(iz)=-Q(z)$, the value $-w_0$ is also regular for $Q$.
Without loss of generality, we can reverse the time and assume that 
$\eps=1$.

We want to describe the fiber $\Psi_t^{-1}(0)$ near $0$, for $t>0$ small.
The assumptions $(i)$, $(ii)$ and $(iii)$ tell us that one can write in a unique way
\begin{equation}
\label{eqndecpsi}
\Psi_t(z)
\; =\;
\Psi_t(0) +Q(z) +R(z)+ t \,r_t(z) 
\end{equation}
where $R$ and $r_t$ are holomorphic maps from $\Om$ to $\m C^m$ 
with
\begin{equation*}
\label{eqnbourrt}
R(z)= o(\|z\|^2)
\;\;{\rm and}\;\;
r_t(0)=0\, .
\end{equation*}
Since $w_0$ is a regular value for $Q$, by Fact \ref{facmulhol}
and \ref{facintqua}, there exist $2^m$ distinct points $v_k\in \m C^m$ indexed by $k=1,\ldots, 2^m$ such that
\begin{equation*}
Q(v_k)\; =\; w_0
\;\;\;{\rm and}\;\;\;
J_Q(v_k)\;\neq\; 0\, ,
\end{equation*}
where $J_Q(z)=\det(DQ(z))\in \m C$ is the Jacobian of the holomorphic map $Q$. We set $M_0:=\max\limits_{k\leq 2^m}\|v_k\|+1$.

Since the map $\Psi$ is singular at the point $(0,0)$, it will be useful to introduce the  auxiliary map 
$\Phi:(-s_0,s_0)\times B_{\m C^m}(0,M_0)\ra \m C^m$ 
given by
\begin{eqnarray*}
\Phi_s(v)&=&
\Psi_{s^2}(s\,v)/s^2
\;\;\;\;\mbox{\rm when}\;\; s\neq 0\\
&=& Q(v)\!-\! w_0
\;\;\mbox{\rm when}\;\; s = 0.
\end{eqnarray*}
When $s_0>0$ is chosen small enough, this map $\Phi$ is well defined 
and  $s\mapsto \Phi_s$ is an analytic family of holomorphic maps.
Indeed, by \eqref{eqndecpsi}, one has 
\begin{eqnarray*}
\Phi_s(v)&=&\Psi_{s^2}(0)/s^2 +Q(v)+ R(sv)/s^2 + r_{s^2}(sv)\;\;\;\mbox{\rm when}\;\; s\neq 0\\
&=& \;\;- w_0\;\;+\;\;Q(v) \;\;\;+\;\;\;\;0\;\;\;\;+\;\;\;\;0
\hspace{7ex}\mbox{\rm when}\;\; s = 0,
\end{eqnarray*}
and each of the four terms is an analytic family of holomorphic maps.

By construction one has 
$$
\Phi_0^{-1}(0)=Q^{-1}(w_0)=\{ v_k\mid k=1,\ldots, 2^m\}.
$$
This map $\Phi$ is particulaly useful since, by the assumption $(iv)$, for all $k\leq 2^m$ 
the differentials $D\Phi_0(v_k)=DQ(v_k)$ are invertible. 
This means that $0$ is a regular value of the map $\Phi_0$. 

By the implicit function theorem, there exists $s'_0>0$, 
such that, for all $k\leq 2^m$,
there exist analytic maps $\ph_k:(-s'_0,s'_0)\ra \m C^m$ 
satisfying
\begin{equation}
\label{eqnphkphi}
\ph_k(0)=v_k
\;\;\;{\rm and}\;\;\;
\Phi_s(\ph_k(s))=0
\;\; \mbox{for all $s\in (-s'_0,s'_0)$}.
\end{equation}

Fix a radius $\rho_0< d(0,\partial \Om)$.
According to Fact \ref{facmulhol}, for $t>0$ small enough
the fiber $\Psi_t^{-1}(0)$ contains at most $2^m$ points in the complex ball $B_{\m C^m}(0,\rho_0)$. By \eqref{eqnphkphi}, these points are exactly the points 
$$
z_k(t):=t^{1/2}\ph_k(t^{1/2}).
$$
For $t$ small enough, the point $z_k(t)$ is real if and only if $v_k$ is real.
By \eqref{eqnfibrea2}, among these $2^m$ points $v_k$,  at most 
$2^{m-1}$ points are real. This proves the bound \eqref{eqnpstpsm}.
\end{proof}

We now come back to the notation of Section \ref{secmulleg}.
We want to construct the unitary deformation $F_t$ of the Fourier transform $F$. 

We introduce the group $U(V)\simeq U(n,\m R)$ of  
unitary transformations of $V\simeq \m C^n$ and its Lie algebra
$$\g u(V):=\{X_0\in {\rm End}(V)\mid \mbox{\rm $X_0$ is antihermitian} \}.
$$
We introduce the euclidean space 
$V_\m R:=\{f\in V\mid \ol f =f\}\simeq \m R^n$, the orthogonal group 
$O(V_\m R):=\{u\in U(V)\mid \ol u=u\}\simeq O(n,\m R)$ and its Lie algebra
$$
\g o(V_\m R)\; :=\;
\{X_0\in \g u(V)\mid \ol{X_0}=X_0\}.
$$ 
We fix an element $X_0\in \g o(V_\m R)$
and we introduce the analytic family 
$$
F_t:=F\circ e^{tX_0}
$$
of unitary deformations of the Fourier transform $F=F_0$.
\vs 

The following lemma introduces an analytic family $\Psi_t$ of holomorphic
maps of $W_0\simeq \m C^{n-1}=\m C^{(p-3)/2}$ which is useful. Here is why.\\ 
$\star$ On the one hand  the fiber $\Psi_t^{-1}(0)$ near $0$ parametrizes 
$T\cap F_t^{-1}T$ near $[\chi_0]$  in the coordinate system \eqref{eqnbetfbe}, for $t$ small.\\ 
$\star$ On the other hand  the family $\Psi_t$ satisfies the assumptions $(i)-(iv)$ of Lemma \ref{lemmuldef}, provided $X_0$ is well chosen.

\bl
\label{lemmulfam} 
$a)$ The following  map $\Psi: \m R\times W_0\ra W_0$ is well defined
\begin{equation*}
\label{eqndefpst}
(t,\be)\mapsto \Psi_t(\be):=\,
e^{tX_0}(\chi_0e^{i\be})\, *\, e^{tX_0}(\chi_0e^{-i\be})\; -\; \chi_0 *\chi_0\, .
\end{equation*}
$b)$ For $t\in \m R$ and $\be \in B_{0,\m R}$, one has the equivalence:
$$\; \chi_0e^{i\be}\in F_t^{-1}T
\Longleftrightarrow 
\Psi_t(\be)=0.$$
$c)$ One can choose $X_0\in \g o(\m R)$ such that the  derivative $w_0:=-\frac{d}{dt}\Psi_t(0)|_{t=0}$ is a regular value for the Hessian $Q:=\tfrac12 D^2\Psi_0(0)$.\\
$d)$ With this choice, the family $\Psi_t$ satisfies the assumptions $(i)\!-\!(iv)$ of Lemma \ref{lemmuldef}. 
In particular there exist a neighborhood $\Om_{\chi_{_0}}$ of $[\chi_0]$ in $\m P(V)$ and a sign $\eps=\pm 1$
such that, for all $t\in (0,1)$ the intersection $T\cap F_{\eps t}^{-1}T$ is transverse in $\Om_{\chi_{_0}}$ and contains at most $ 2^{n-2}$ points in $\Om_{\chi_{_0}}$.
\el

\begin{proof}
$a)$ The proof is as in Lemma \ref{lemmulmul}.$a$. 
We want to check that, for $\be$ in $W_0$, the function $\Psi_t(\be)$ also belongs to $W_0$.
Since the functions $f_{\pm \be}=\chi_0 e^{\pm i\be}$ are odd,
the functions $e^{tX_0}f_{\pm \be}$ are also odd
and the function $\Psi_t(\be)$
is even.\\
Moreover, since the unitary map $e^{tX_0}$ is orthogonal, one also has
$$
\Psi_t(\be)(0)\;=\; -\|e^{tX_0}f_{\be}\|_{\ell^2}^2+ \|\chi_0\|_{\ell^2}^2
\;=\; -\|f_{\be}\|_{\ell^2}^2+ \|\chi_0\|_{\ell^2}^2\;=\;0
$$
Since  the function $e^{tX_0}f_{\be}$ and its Fourier transform 
$\widehat{e^{tX_0}f_{\!\be}}$ are odd, one has
$$\textstyle
\tfrac{1}{\sqrt{p}}\,\widehat{\Psi_t(\be)}(0)\;=\; \widehat{e^{tX_0}f_{\be}}(0)\; \widehat{e^{tX_0}f_{\!-\be}}(0)\;-\; (\widehat{\chi_0}(0))^2\;=\;0
$$
This proves that the function $\Psi_t(\be)$ belongs to $W_0$.

$b)$ The proof is as in Lemma \ref{lemmulmul}.$b$.  

$c)$ This follows from the following four assertions:\\ 
The derivative $w_0$ is equal to $w_0=-2\,X_0(\chi_0)*\chi_0$.\\ 
The map $\g o(V_\m R)\ra V_{0,\m R};X_0\mapsto \al:= X_0(\chi_0)$ is onto.\\
The map $V_{0,\m R}\ra W_{0,\m R} ;\al\mapsto \be:=\al *\chi_0$ is an isomorphism.\\
The real vector space $W_{0,\m R}$ is not contained in the set $Z\subset W_0$ 
of singular values of the holomorphic quadratic map $Q:W_0\ra W_0$.

$d)$ The first statement follows from Point $c)$ and from Lemma \ref{lemasupsi}. 
The last statement follows from Point $b)$ and from Lemma \ref{lemmuldef}.
\end{proof}

This ends the proof of Proposition \ref{protrach0} and hence this also ends 
the proof of Theorem \ref{thmexicoh}.

\section{Counting $\mc C$-functions}
\label{seccoucoh}

The aim of this Chapter is to prove 
Theorem \ref{thmcouequ}
which is an extension of Theorem \ref{thmcouodd}.
This extension gives a formula for the number of $(H,c)$-equivariant $\mc C$-functions.

In Section \ref{secboucoh}, we give a short and new proof of the finiteness of the set of $\mc C$-functions $f$ with $f(1)=1$.

In Section \ref{seccouequ}, we explain the deformation argument
that allows one to compute the number of $(H,c)$-equivariant $\mc C$-functions.

\subsection{Finiteness of the set of $\mc C$-functions}
\label{secboucoh}

The following proposition is due to Biro.

\bp
\label{proboucoh}  {\bf (\cite{BiroJNT})}
Let $p$ be a prime number. Then the set of $\mc C$-functions $f$ on $\m F_p$ such that $f(1)=1$ is finite.
\ep

We will give a different proof of this proposition 
whose arguments will be reused in the sequel. 
We first recall a lemma due to Chebotarev that we call the little Chebotarev theorem.
 
\bl
\label{lemlitche} {\bf (\cite{LenstraStevenhagen})}
Let $p$ be a prime number and $\zeta_p:=e^{2i\pi/p}$. For any subset $A$, $B$ of $\m F_p$ with same cardinality,
the matrix $(\zeta_p^{jk})_{ \substack{j\in A\\k\in B}}$ is invertible
\el

The little theorem of Chebotarev has been reinterpreted by Biro and Tao
as the following equivalent statement.

\bl
\label{lemtaocyc} {\bf (see \cite[p.2]{BiroLev} and  \cite{TaoCyclic})}
Let $p$ be a prime number. For every non zero 
function $f$ on $\m F_p$, one has 
\begin{equation}
\label{eqntaoceb}
\#\,{\rm supp}(f)+\#\,{\rm supp}(\widehat{f})
\;\geq\; 
p\!+\!1\, .
\end{equation}
\el

This lemma is the key ingredient towards the following proposition due to Haagerup. Let $E=\m C^{\m F_p}$ be the set of functions on $\m F_p$ endowed with the sup-norm $\|.\|_\infty$
and $E_{\geq 1}$ be the complementary of the open ball of radius $1$, 
$$
E_{\geq 1}:=\{f\in E\mid \|f\|_\infty \geq 1\}.
$$

\bp 
\label{prohaaphi} {\bf (\cite{HaagerupFourier})}
Let $p$ be a prime number. 
Then the map 
\begin{eqnarray*}
	\Phi: E_{\geq 1}\times E_{\geq 1}
	&\longrightarrow &
	E\times E\\
	(f,g)
	&\mapsto& (fg,\widehat{f}\, \widehat{\widecheck{g}})
\end{eqnarray*}
is a proper map.
\ep

``Proper'' means that the inverse image of a compact set is compact.

\br 
The same would be true with the map 
$\Psi:(f,g)\mapsto (fg,\widehat{f}\, \widehat{g})$.
But the map $\Phi$ wich involves $\widecheck{g}$ 
is more useful for us because of 
Lemma \ref{lemfoucohn}.
\er

\begin{proof}[Proof of Proposition \ref{prohaaphi}]
	Assume by contradiction that there exists 
	sequen\-ces $f_n$, $g_n$ in $E_{\geq 1}$ with 
${\rm max}(\| f_n\|_\infty,\| g_n\|_\infty)$ going to $\infty$
such that $\Phi(f_n,g_n)$ converges. After extracting a subsequence, the renormalized functions converge
$$
u_n:=\frac{f_n}{\| f_n\|_\infty}\ra u_\infty
\;\; {\rm and} \;\;
v_n:=\frac{g_n}{\| g_n\|_\infty}\ra v_\infty\, ,
$$
and the limit functions $u_\infty$ and $v_\infty$ are non zero functions. Since the product 
$\|f_n\|_\infty\|g_n\|_\infty$ goes to $\infty$,
these limits satisfy
$$
u_\infty v_\infty=0
\;\; {\rm and} \;\;
\widehat{u}_\infty\widehat{\widecheck{v}}_\infty=0.
$$
In particular, one has the inequality
$$
\#\,{\rm supp}(u_\infty)+\#\,{\rm supp}(v_\infty)+\#\,{\rm supp}(\widehat{u}_\infty)+\#\,{\rm supp}(\widehat{v}_\infty)
\;\leq\; 
2p ,
$$
which contradicts Lemma \ref{lemtaocyc}.
\end{proof}

\begin{proof}[Proof of Proposition \ref{proboucoh}]
We introduce the complex affine subspace of dimension $p\!-\! 2$
in $E:=\m C^{\m F_p}$,
$$
E_1:=\{f\in \m C^{\m F_p]}\mid f(0)=0\; {\rm and}\; f(1)=1\},
$$
We have seen in Proposition \ref{prohaaphi} that the map
\begin{eqnarray}
\label{eqnphifxf}
\Phi: E_1\times E_1	&\longrightarrow &
	E\times E\\
	(f,g)
	&\mapsto& (fg,\widehat{f}\, \widehat{\widecheck{g}})
\end{eqnarray}
is a proper map. 
By Lemma \ref{lemfoucohn}, 
a function $f$ in $F_1$ is a $\mc C$-function if and only if there exists a (unique) function $g$ in $E_1$ such that
$\Phi(f,g)=({\bf 1}_{\m F_p^*},{\bf 1}_{\m F_p^*})$.
By Proposition \ref{prohaaphi} the complex algebraic subvariety $\Phi^{-1}({\bf 1}_{\m F_p^*},{\bf 1}_{\m F_p^*})$ of 
$E_1\times E_1$ is compact. 
Therefore it is finite.
\end{proof}

\subsection{Counting equivariant $\mc C$-functions}
\label{seccouequ}

 We want to prove the following.

\bt 
\label{thmcouequ}
Let $p$ be an odd prime, $H\subset
\m F_p^*$ be a subgroup of index $n$, 
and $c:H\ra \m C^*$ be a non-trivial character. 
Then, counted with multiplicities,
the number of $(H,c)$-equivariant  $\mc C$-functions $f$ on $\m F_p$ with $f(1)=1$, is equal to the binomial coefficient $\binom{2n-2}{n-1}$.
\et

\begin{proof}[Proof of Theorem \ref{thmcouodd}]
Just apply Theorem \ref{thmcouequ} with the group $H=\{\pm 1\}$ and with the  unique odd character $c$ of $H$. The $(H,c)$-equivariant functions $f$  are nothing but the odd $\mc C$-functions and the index $n$ is $n=(p\!-\! 1)/2$. 
\end{proof}

The key point will be a deformation argument.
For this argument
one needs a defining map $\Phi$ as in \eqref{eqnphifxf} which is not only proper but also dominant so that the number of points in its fibers, counted with multiplicities is constant, 
and so that we can perform the counting for a simpler fiber.

We introduce the vector space $V=V_{H,c}$ of $(H,c)$-equivariant functions $f$ on $\m F_p$ as in \eqref{eqnequfun},
and its affine subspace 
\begin{equation*}
\label{eqnv1fvhc}
V_1:=\{f\in V_{H,c}\mid f(1)=1\}.
\end{equation*}
Since $c$ is non trivial, all these $(H,c)$-equivariant functions  vanish at $0$. 
The conjugate affine space is 
$$
\ol{V}_1:=\{g\in V_{H,\ol{c}}\mid g(1)=1\}.
$$
We also introduce the space $V_0$ of $H$-invariant functions $f_0$ on $\m F_p$ such that $f_0(0)=0$, and the affine space
\begin{equation*}
\label{eqnwphpsi}
W_1:=\{(f_0,g_0)\in V_0\times V_0\mid f_0(1)=1
\;\;{\rm and}\;\;
\textstyle\sum\limits_{x\in \m F_p} f_0(x)=
\sum\limits_{x\in \m F_p}g_0(x)\;\}.
\end{equation*}

\bl 
\label{lemcouequ}
The map
\begin{eqnarray}
\label{eqnphv1v1}
\Phi: V_1\times {\ol V}_1	&\longrightarrow &
W_1\\
\nonumber
(f,g)
&\mapsto& (fg,\widehat{f}\, \widehat{\widecheck{g}})
\end{eqnarray}
is a well-defined dominant and proper map. 
\el

Note that $\dim(W_1)=2\dim(V_1)=2n-2$.

\begin{proof}[Proof of Lemma \ref{lemcouequ}]
We first check that $\Ph(V_1\times \ol{V}_1)$ is included in $W_1$.  Let $f$ be in $V_1$, let $g$ be in $\ol{V}_1$.
Since $f$ is $(H,c)$ equivariant and $g$ is $(H,\ol{c})$-equivariant, the function $fg$ is $H$-invariant
and vanish at $0$.
Similarly, since $\widehat{f}$ is $(H,\ol{c})$ equivariant and $\widehat{\widecheck{g}}$ is $(H,c)$-equivariant, the function $\widehat{f}\, \widehat{\widecheck{g}}$ is $H$-invariant
and vanish at $0$.
Since $f(1)=g(1)=1$, one also has $fg(1)=1$
The last equality $\sum\limits_{x\in \m F_p} fg(x)=
\sum\limits_{x\in \m F_p}\widehat{f}\, \widehat{\widecheck{g}}(x)$ is the  Parseval formula
for finite Fourier transform.

The fact that the map $\Phi$ is proper follows from Proposition \ref{prohaaphi}.

The fact that $\Phi$ is dominant follows from properness and 
equidimensionality. One can also exhibit an explicit point 
$(f,g)$
at which the differential $D\Phi(f,g)$ is an isomorphism. 
For instance the points in  Lemma \ref{lemcoufib}.$e$.
\end{proof}

In the following lemma we study in great detail the fiber of $\Phi$ over the 
point $({\bf 1}_H,{\bf 1}_H)$.

\bl 
\label{lemcoufib} Let $H$ be a subgroup of 
$G:=\m F_p^*$ and $c$ be a non trivial character of $H$.
Let $\Phi$ be the map \eqref{eqnphv1v1}.\\ 
$a)$ For every $H$-invariant subsets $A$, $B$ of $\m F_p^*$
such that 
\begin{equation}
\# A+\# B=\# G+\# H,
\end{equation}
there exists a unique function $f_{A,B}\in V_1$ with support $A$
whose Fourier transform $\widehat{f}_{A,B}\in V_{H,\ol{c}}$ has support $B$.\\
$b)$ Assume now that both $A$ and $B$ contain $H$. We set $g_{A ,B}\in\ol{V}_1$ to be the function 
$g_{A ,B}:=\ol{f}_{A',-B'}$ where $A'$ and $B'$ are defined by $A\cup A'=B\cup B'=G$ and $A\cap A'=B\cap B'=H$.
Then one has 
$$
\Phi(f_{A,B},g_{A,B})=({\bf 1}_H,{\bf 1}_H).
$$
$c)$ Conversely every point in the fiber
 $\Phi^{-1}({\bf 1}_H,{\bf 1}_H)$ is one of these points 
$(f_{A,B},g_{A,B})$.\\
$d)$ The number of points in this fiber $\Phi^{-1}({\bf 1}_H,{\bf 1}_H)$ is equal to $\binom{2n-2}{n-1}$.\\
$e)$ The map $\Phi$ is non-degenerate at each of the points $(f_{A,B},g_{A,B})$ of the fiber $\Phi^{-1}({\bf 1}_H,{\bf 1}_H)$.
\el

\begin{proof}[Proof of Lemma \ref{lemcoufib}]
Let $d=p\! -\! 1$, and write $d=n \,d_H$ where $n$ is the index of $H$
and $d_H$ the order of $H$. 

$a)$ For a $H$-invariant subset
$A$ of $G$ of cardinality $n_Ad_H$, we denote by   $\m C[A]_{H,c}$ the vector space of $(H,c)$-equivariant functions
on $A$.  It has dimension $n_A$.

As a consequence of  Inequality \eqref{eqntaoceb}, for any $H$-invariant subsets $A_0$ and $B_0$ of $G$ with $n_{A_0}+n_{B_0}=n$,
the map 
$$\m C[A_0]_{H,c}\ra \m C[B_0^c]_{H,\ol{c}}: f\ra \widehat{f}|_{B_0^c}$$ 
is an isomorphism.

Therefore, for any $H$-invariant subsets $A$ and $B$ of $G$ with $n_{A}+n_{B}=n+1$,
the map 
$$\m C[A]_{H,c}\ra \m C[B^c]_{H,\ol{c}}: f\ra \widehat{f}|_{B^c}$$ 
has a one dimensional kernel. Moreover, a non-zero function $f_{A,B}$ in the kernel does not vanish on $A$; one can normalize it so that 
$f_{A,B}(1)=1$. Similarly its Fourier transform $\widehat{f}_{A,B}$
does not vanish on $B$.

$b)$ Let $f:=f_{A,B}$ and $g:=g_{A,B}$. By construction the product $fg$ is $H$-invariant, is supported by $H$ and satisfies
$fg(1)=1$. Therefore one has $fg={\bf 1}_H$.

Similarly the product $\widehat{f}\widehat{\widecheck{g}}$ 
is $H$-invariant and is supported by $H$. Therefore 
one has $\widehat{f}\widehat{\widecheck{g}}=\la {\bf 1}_H$ 
for some constant $\la$. This constant $\la$  is equal to $1$ 
by the Parseval-Plancherel formula.

$c)$ Conversely, let $(f,g)\in V_1\times\ol{V}_1$ such that 
$fg={\bf 1}_H$ and $\widehat{f}\widehat{\widecheck{g}}= {\bf 1}_H$.
Set $A:={\rm supp}(f)$, $B:={\rm supp}(\widehat{f})$,  
$A':={\rm supp}(g)$ and $B':={\rm supp}(\widehat{\widecheck{g}})$.
By assumption, one has $A\cap A'=B\cap B'=H$.
In particular, one has 
$$
n_A+n_{A'}\leq  n+1
\;\;{\rm and}\;\;
n_{B}+n_{B'}\leq  n+1.
$$
By inequality \eqref{eqntaoceb}, one has 
$$
n_A+n_{B}\geq  n+1
\;\;{\rm and}\;\;
n_{A'}+n_{B'}\geq  n+1.
$$ 
Therefore all these inequalities are equalities and hence,
by point $a)$, one has $f=f_{A,B}$ and $g=\ol{f}_{A',-B'}=g_{A,B}$.

$d)$ By Point $c)$, the fiber 
$\Phi^{-1}({\bf 1}_H,{\bf 1}_H)$ is in bijection with the set of 
pairs $(A\smallsetminus H,B\smallsetminus H)$ of $H$-invariant subsets of $G\smallsetminus H$
such that $n_{A\smallsetminus H}+n_{B\smallsetminus H}=n-1$. Their total number is $\binom{2n-2}{n-1}$ as announced.

$e)$ Fix a point $(f,g)=(f_{A,B},{\ol{f}_{A',-B'}})$ in the fiber
$\Phi^{-1}({\bf 1}_H,{\bf 1}_H)$. We want to prove that the differential $D\Phi(f,g)$ is injective.
The tangent space of the source is the space of 
couples $(\ph,\psi)$ of $(H,c)$-equivariant functions 
such that $\ph(1)=\psi(1)=0$. 
Assume that $(\ph,\psi)$ is in the kernel of $D\Phi(f,g)$.
The formula for the differential is 
$$
D\Phi(f,g)(\ph,\psi)=(f\psi+g\ph, \widehat{f}\widehat{\widecheck{\psi}}+ 
\widehat{\widecheck{g}}\widehat{\ph})=0.
$$
Since the functions $f\psi$ is supported by $A\smallsetminus H$
and the function $g\ph$ is supported by $A'\smallsetminus H$,
one gets $f\psi=g\ph=0$.
Since $f$ does not vanish on $A$ and $g$ does not vanish on $A'$, this proves that 
$$
{\rm supp}(\ph)\subset A\smallsetminus H
\;\;{\rm and}\;\;
{\rm supp}(\psi)\subset A'\smallsetminus H.
$$
A similar argument proves that 
$\widehat{f}\widehat{\widecheck{\psi}}=
\widehat{\widecheck{g}}\widehat{\ph}=0$ and that
$$
{\rm supp}(\widehat{\ph})\subset B
\;\;{\rm and}\;\;
{\rm supp}(\widehat{\psi})\subset -B'.
$$
In particular one gets 
\begin{eqnarray*}
\#{\rm supp}(\ph)+\#{\rm supp}(\widehat{\ph}) &\leq& \# G,\\
\#{\rm supp}(\psi)+\#{\rm supp}(\widehat{\psi}) &\leq& \# G.
\end{eqnarray*}
Therefore, by Lemma \ref{lemtaocyc}, one has $\ph=\psi=0$.

This proves that the differential $D\Phi(f,g)$  is an isomorhism.
\end{proof}

\begin{proof}[Proof of Theorem \ref{thmcouequ}]
	Since the map $\Phi$ in \eqref{eqnphv1v1} is proper and dominant,
	the number of points in its fibers $\Phi^{-1}(w)$, counted with multiplicities is constant. See Fact \ref{facmulhol}. See also \cite[Sec. 4]{HaagerupFourier}. 
	
	Remembering Lemma \ref{lemfoucohn} we want to prove that, counted with multiplicity,
	the number of points in the fiber
	$\Phi^{-1}({\bf 1}_{\m F_p^*},{\bf 1}_{\m F_p^*})$
	is equal to  $\binom{2n-2}{n-1}$.
	
	It is equivalent to prove it for the fiber
	$\Phi^{-1}({\bf 1}_H,{\bf 1}_H)$. This was done in the previous Lemma \ref{lemcoufib}.
\end{proof}

\br
It would be nice to have a similar counting formula for all $\mc C$-functions $f$ on $\m F_p$ with $f(1)=1$. 
\er

\section{Examples}
\label{secexalis}

We report in this chapter a few examples, relying on numerical experiments, that emphasize on the one hand the arithmetic complexity 
of these new unimodular  $\mc C$-functions that we have proven to exist in this article, and on the other hand some kind of ``uniqueness''
of these new unimodular $\mc C$-functions.

\subsection{When $d=7$}
\label{secexad07}

There are exactly
$\binom{4}{2}=6=1\!\times\!4+\!2$
odd $\mc C$-functions with $f(1)=1$
on $\m F_7$ counted with mutiplicities. All of them being Dirichlet characters.
All of them except the Legendre character $\chi_0$ having multiplicity $1$. 

All of them are unimodular.
As expected, this number is greater than  $2^2=4$ which is the lower bound predicted by Theorem \ref{thmchobis}.

\subsection{When $d=9$}
\label{secexad09}

There are exactly
$18$ odd $\mc C$-functions with $f(1)=1$
on the cyclic group $C_9$. All of them having multiplicity one.
Among those functions, there are no Dirichlet characters.
These $18$ functions are Galois conjugate.

The total number of unimodular odd $\mc C$-functions
is 
$12$.
As expected, this number is greater than  $2^3=8$ which is the lower bound predicted by Theorem \ref{thmchobis}.

\subsection{When $d=11$}
\label{secexad11}

There are exactly
$\binom{8}{4}=70=1\!\times\!16+\!4\!+\!5\!+\!5\!+\!40$ odd $\mc C$-functions on $\m F_p$, counted with multiplicities. All of them except one having multiplicity one. The exception is the  Legendre character $\chi_0$ which has multiplicity $16$.
Among those functions, there are also the $4$ Dirichlet characters of order $10$.
There are also the $5$ real-valued functions which are Galois conjugate to
the function $k\mapsto (\!\frac{k}{11}\!)\,(c(k^2)+2\,c(4k^2))$ where 
$(\!\frac{k}{11}\!)$ is the Legendre symbol and where $c(m):={\rm cos}(2m\pi/11)$. There are also the $5$ inverses of these functions.

The last $40$ functions are Galois conjugate. 

The total number of unimodular odd $\mc C$-functions, counted with multiplicities,
is 
$30=1\!\times\! 16+\!4\!+\!0\!+\!0\!+\!10$.
As expected, this number $30$ is greater than  $2^4=16$ which is the lower bound predicted by Theorem \ref{thmchobis}.

\subsection{When $d=13$}
\label{secexad13}

There are exactly
$\binom{10}{5}=252=2\!+\!4\!+\!6\!+\!96\!+\!144$ odd $\mc C$-functions on 
$\m F_p$ all of them with multiplicity $1$.
Among those, are the $2$ Dirichlet characters
of order $4$ and the $4$ Dirichlet characters of order $12$.
Among those are also $6$ other odd functions $f$ such that, for some $\eps=\pm 1$, one has $f(5x)=i^{\eps} f(x)$ for all $x$.
The remaining functions form two Galois conjugacy classes,
the first one contains $96$ functions, the last one contains $144$ functions. 

The total number of unimodular odd $\mc C$-functions
is 
$60=2\!+\!4\!+\!6\!+\!0\!+\!48$.
As expected, this number $60$ is greater than  $2^5=32$ which is the lower bound predicted by Theorem \ref{thmchobis}.

\appendix 

\section{Appendix: Biunimodular functions}
\label{secbiufun} 

Harvey Cohn's problem was developped simultaneously with the very popular problem which is called the 
{\it biunimodular functions problem}.
The aim of this appendix is to recall part of the history of this problem and to explain in Theorem \ref{thmexibiu} new constructions of biunimodular 
functions that can be obtained 
by the method we have developped in this paper.

\subsection{Classical biunimodular functions} 
\label{secpreres}

Let $d$ be an integer.  
A complex valued function $f$ on the cyclic group $C_d:=\m Z/d\m Z$ is called {\it biunimodular} if  $|f(\ell)|=1$ for all $\ell\in C_d$. and if
\begin{equation} 
	\label{eqnpcyclic0}\textstyle
	\sum_{k\in \m Z/d\m Z} f(k\!-\!\ell)\overline{f(k)}\; =\; 0
	\;\;\;\;\mbox{for all $\ell\neq 0$.}
\end{equation}
Geometrically, Condition \eqref{eqnpcyclic0} means that the translates of $f$ form an
orthogonal basis of $\ell^2(C_d)$.
It is equivalent to require
\begin{equation} 
	\label{eqnpcyclic1}\textstyle
	|\widehat{f}(\ell)|=1	\;\;\;\;\mbox{for all $\ell\in C_d$.}
\end{equation}
This notion was introduced by
Per Enflo in the $80$'s in relation with the ``circulant complex Hadamard matrices''.
The simplest examples of biunimodular functions are the gaussian functions. 
We recall that, when $d$ is odd, a  gaussian function  on the cyclic group $C_d=\m Z/d\m Z$ is a function of the form 
\begin{equation}
	\label{eqngaufun}
	g_{m,a}:x\mapsto e^{2i\pi m(x-a)^2/d}
\end{equation}
for some $m$, $a$ in $\m Z/d\m Z$, $m$ coprime to $d$.

Bj\"{o}rck and Saffari
classified in 1995 in \cite{BjorckSaffari}, when $d=p>2$ is prime,
the biunimodular functions that are $(\m F_p^*)^2$-invariant. They found:

When $p\equiv 1$ mod $4$ there are, up to proportionality, two such functions 
\begin{equation*}
	\label{eqnphp1m4}
	h_\eps =\de_0+ {\rm cos}\theta \,{\bf 1}_{\m F_p^*}+ i\,\eps\,{\rm sin}\theta\, \chi_0,
\end{equation*}
where $\chi_0$ is the Legendre character, 
$ {\rm cos}\theta :=\frac{1}{\sqrt{p}+1}$ and $\eps=\pm 1$.

When $p\equiv 3$ mod $4$ there are,  up to proportionality, four such functions 
\begin{equation*}
	\label{eqnphp3m4}
	h_{\eps} =e^{i\,\eps_1 \theta}\de_0+ cos\theta\, {\bf 1}_{\m F_p^*}+
	i\,\eps_2\,{\rm sin}\theta \,\chi_0,
\end{equation*}
where $\chi_0$ is the quadratic Dirichlet character, 
$ {\rm tan}\theta =\sqrt{p}$ and $\eps_1,\eps_2=\pm 1$.

We introduce the translates of these functions where $a\in \m F_p$
\begin{equation} 
	\label{eqnph13m4}
	h_{\eps,a}:x \mapsto h_\eps(x -a).
\end{equation} 
For all prime $p\geq 7$, this gives rise to either $2p$ or $4p$  new biunimodular functions  that we call the Bj\"{o}rck-Saffari functions.

Haagerup proved in 2008 in \cite{HaagerupFourier} that, 
when $p$ is prime, the set of biunimodular functions with $f(0)=1$ is finite.

When $p\equiv 1$ mod $3$ is prime, biunimodular functions 
on $\m F_p$ which are $(\m F_p^*)^3$-invariant
have been classified by Gabidulin and Shorin
in \cite{GabidulinShorin} and  by Bj\"{o}rck and Haagerup
in \cite{BjorckHaagerup}. This gave rise to new biunimodular functions for all prime $p\equiv 1$ mod $3$, $p\geq 13$. 
Their method relies on explicit computations 
that are possible since the $(\m F_p^*)^3$-invariant functions $f$ with $f(0)=1$ depend only on three parameters.

Using the same strategy when  $p\equiv 1$ mod $5$, a few computer assisted constructions of $(\m F_p^*)^5$-invariant biunimodular functions were given for a few values of $p$
in \cite{LoidreauShorin}: $p=31$ and $p=61$.

Biunimodular functions are studied under various names like CAZAC, PSK perfect sequences, polyphase sequences with optimum correlation and so on. See \cite{BenedettoCordwellMagsino}, \cite{FuhrRzeszotnik} and \cite{Nicoara} 
for strongly related topics. 

\subsection{Constructing biunimodular functions} 
\label{secnewres}

The method we developped  in this paper allows us to have a new insight 
in this problem. For instance we will prove: 

\bt 
\label{thmexibiu}
Let  $p\geq 11$ be prime. There exist biunimodular functions on 
$\m F_p$ which are proportional neither to gaussian
nor to Bj\"{o}rck-Saffari functions.
\et

The proof will be similar to the proof of Theorem \ref{thmexicoh}.
Let  $V:=\ell^2(\m F_p)$ be the $p$-dimensional Hilbert space of functions $f$ on $\m F_p$.
Let 
$
\m P(V)\simeq \m C\m P^{p-1}
$ 
be the projective space  of $V$, let $T$ be the Clifford torus 
\begin{equation}
	\label{eqnclitor2}
	T:=\{[f]\in \m P(V)\mid \;|f(x)|=|f(0)|\; \mbox{\rm for all $x\in \m F_p$}\}\simeq \m T^{p-1}.
\end{equation}
and $F:f\mapsto \widehat{f}$ be the Fourier transform
given by \eqref{eqnfoufor}.

\begin{proof}[Proof of Theorem \ref{thmexibiu}]
	We will apply Theorem \ref{thmchobis} to the unitary transformation $F$ and the torus $T$.
	The $(p-1)p$ gaussian functions $g_{m,a}$ 
	as in \eqref{eqngaufun} belong to the intersection $T\cap F^{-1}T$. 
	The $2p$ or $4p$ functions 
	$h_{\eps,a}$ as in \eqref{eqnph13m4} also 
	belong to the intersection $T\cap F^{-1}T$.
	
	Assume, by contradiction that the intersection $T\cap F^{-1}T$ contains only gaussian and Bj\"{o}rck-Saffari functions.
	Then, by the following Propositions \ref{protragau} and \ref{protrabjo}, this intersection is transverse. Therefore 
	Theorem \ref{thmchobis} predicts the existence of 
	at least $2^{p-1}$ intersection points counted with multiplicity.
	Since $p\geq 11$, one has $2^{p-1}>(p-1)p+4p$, there must exist another intersection point. This is the contradiction we are looking for.
\end{proof}

The proof  relied on  transversality properties that will be proven in Sections \ref{sectragau} and \ref{sectrabjo}.

\subsection{Transversality at gaussian functions}
\label{sectragau}

Let $p\geq 3$ be prime. 
The aim of this section is to prove the transversality of the intersection $T\cap F^{-1}T$ at the gaussian functions.
To simplify notation we will only deal with one gaussian function.

\bp 
\label{protragau}
Let $p\geq 3$ be prime and $g_0$ be the gaussian  function on $\m F_p$,
$x \mapsto g_0(x ):=e^{2i\pi x ^2/p}$.
Then the intersection $T\cap F^{-1}T$ is transverse at $[g_0]$.
\ep

\begin{proof}[Proof of Proposition \ref{protragau}]
	We compute the Fourier transform: for  $ x $ in $\m F_p$, 
	\begin{equation*}
		\widehat{g}_0( x )
		\;=\;\eps _p\; \overline{g_0}( x /2),
	\end{equation*}
	where $\eps_p:=\tfrac{G(\chi_0)}{\sqrt{p}}$ and 
	where $1/2$ is the inverse of $2$ in $\m F_p$.
	
	Note that $\eps_p=1$ when $p\equiv 1$ mod $4$ and 
	$\eps_p=i$ when $p\equiv 3$ mod $4$.  
	\vs 
	
	$\star$ {\bf We describe the tangent space to the projective space.}
	
	It will be convenient to use a parametrization of a neighborhood of $[g_0]$ 
	in $\m P(V)$ by the vector space 
	$V_o:=\{\ph\in \m C^{\m F_p}\mid \ph(0)=0\}$ given by 
	\begin{equation*}
		\label{eqncooafa1}
		{\ph}\mapsto [g_{\ph}]
		\;\;{\rm where}\;\; 
		g_{\ph}=
		\textstyle\left({\bf 1}_{\m F_p}\!+\!\ph\right)g_0\, .
	\end{equation*}
	This gives an identification of $V_o$
	with the tangent space of $\m P(V)$ at the point $[g_0]$, thanks to the formula
	\begin{equation*}
		\label{eqntanpv1}
		{\ph}\mapsto v_{\ph}:=\frac{d}{d\eps}[g_{\eps{\ph}}]|_{\eps=0}
		\;\in\; T_{[g_0]}\m P(V) .
	\end{equation*}
	
	$\star$ {\bf We describe the tangent space to the torus $T$}.
	Since $\ph(0)=0$, the linear condition defining the tangent space of $T$ at the point $[g_0]$ is
	\begin{equation}
		\label{eqntanga2}
		\textstyle
		{\rm Re}(\ph) \;=\; 0.
	\end{equation}

	$\star$ {\bf We describe the tangent space to the torus $F^{-1} T$}.
	Using \eqref{eqnchieig}, one computes in our coordinate system
	\begin{eqnarray*}
		\widehat{g}_{\ph}
		&=&
		\textstyle
		\left({\bf 1}_{\m F_p}+U\ph\right)\widehat{g}_{0} 
		\;\;{\rm where}\\
		U\ph(x)
		&=&\textstyle
		\eps_p^{-1} \, g_0( x/2)\sum_{y\in\m F_p}e^{2i\pi xy/p}e^{2i\pi y^2/p}\ph(y)\\
		&=&\textstyle 
		\eps_p^{-1} \, \sum_{y\in\m F_p}g_0(x/2+y)\,\ph(y)\\
		&=& \eps_p^{-1} \,(g_{0} *\ph)(-x/2),
		\;\; \mbox{\rm for all $x$ in $\m F_p$.}
	\end{eqnarray*}
	The linear condition defining the tangent space of $F^{-1}T$ at the point $[g_0]$ is
	\begin{equation}
		\label{eqntanga3}
		\textstyle
		{\rm Re}(U\ph) \;\; 
		\mbox{\rm is constant on $\m F_p$}.
	\end{equation}	
	$\star$ {\bf We check the transversality of these tangent spaces}.
	We want to prove that a function $\ph\in V_o$ belonging to both tangent spaces is zero. 
	By \eqref{eqntanga2} one can write $\ph=i\psi$ with $\psi$ real valued. We set $g_0=\al_0+i\be_0$
	where
	$$
	\al_0(x)=\cos(2\pi x^2/p)
	\;\; {\rm and}\;\;
	\be_0(x):=\sin(2\pi x^2/p).
	$$
	and we distinguish two cases.
	\vs 
	
	{\bf First case: when  $p\equiv 1$ mod $4$}. Equation \eqref{eqntanga3} can be rewritten as 
	\begin{equation*}
		\label{eqntanga3a}
		\textstyle
		\be_0 *\psi \;\; 
		\mbox{\rm is constant on $\m F_p$}.
	\end{equation*}	
	or equivalently 
	\begin{equation*}
		\label{eqntanga4a}
		\textstyle
		\widehat{\be}_0 \widehat{\psi} \;\; 
		\mbox{\rm is zero on $\m F_p^*$}.
	\end{equation*}
	Since the function $\widehat{\be}_0(x)=-\be_0(x/2)$ does not vanish on $\m F_p^*$, this implies that $\widehat{\psi}$ is zero on $\m F_p^*$.
	Therefore, since $\sum_y\widehat{\psi}(y)=\sqrt{p} \,\psi(0)=0$,
	one gets $\widehat{\psi}=0$ and $\psi=0$, as required.
	\vs 
	
	{\bf Second case: when  $p\equiv 3$ mod $4$}. Equation \eqref{eqntanga3} can be rewritten as 
	\begin{equation*}
		\label{eqntanga3b}
		\textstyle
		\al_0 *\psi \;\; 
		\mbox{\rm is constant on $\m F_p$}.
	\end{equation*}	
	or equivalently 
	\begin{equation*}
		\label{eqntanga4b}
		\textstyle
		\widehat{\al}_0 \widehat{\psi} \;\; 
		\mbox{\rm is zero on $\m F_p^*$}.
	\end{equation*}
	Since the function $\widehat{\al}_0(x)=i\,\al_0(x/2)$ does not vanish on $\m F_p^*$, this implies that $\widehat{\psi}$ is zero on $\m F_p^*$.
	Therefore, as above, the function $\psi$ is $0$.
\end{proof}

\subsection{Transversality at Bj\"{o}rck functions}
\label{sectrabjo}

Let $p\geq 3$ be prime. 
The aim of this section is to prove the transversality of the intersection $T\cap F^{-1}T$ at the Bj\"{o}rck-Saffari functions.
To simplify we will only deal with one such function and choose notation that works both for $p\equiv  \pm 1$ mod $4$.

When $p\equiv 1$
mod $4$, we set  $\theta_0=0$ and 
$\theta= \arccos(\frac{1}{\sqrt{p}+1})$.

When $p\equiv 3$
mod $4$, we set  $\theta_0=\theta= \arctan(\sqrt{p})$.

\bp 
\label{protrabjo}
Let $h_0$ be the Bj\"{o}rck-Saffari function on $\m F_p$,
\begin{equation}
	\label{eqntrabjo}
	h_0 \;=\;
	e^{i\,\theta_0}\de_0+ \cos\theta\, {\bf 1}_{\m F_p^*}+
	i\,{\rm sin}\theta \,\chi_0.
\end{equation}
Then the intersection $T\cap F^{-1}T$ is transverse at $[h_0]$.
\ep

\begin{proof}[Proof of Proposition \ref{protrabjo}]
	We compute the Fourier transform:  
	\begin{equation}
		\label{eqntrabjo2}
		\eps_p^{-1}\,\widehat{h}_0
		\;=\;e^{-i\,\theta_0}\de_0+ \cos\theta\, {\bf 1}_{\m F_p^*}+
		i\,{\rm sin}\theta \,\chi_0,
	\end{equation}
	where  $\eps_p=1$ when $p\equiv 1$ mod $4$ and 
	$\eps_p=i$ when $p\equiv 3$ mod $4$.
	\vs 
	
	$\star$ {\bf We describe the tangent space to the projective space.}
	The set  $B$ of Dirichlet characters on $\m F_p$ is a  basis of $V_{o}$. 
	It will be convenient to use the following coordinates system 
	${\bf a}=({\bf a}_\chi)_{\chi\in B}$ of 
	$\m P(V)$ in the neighborhood of $[h_0]$ 
	where the coordinates ${\bf a}_\chi$ are complex numbers. It is given by 
	\begin{equation*}
		\label{eqncooafa2}
		{\bf a}\mapsto [h_{\bf a}]
		\;\;{\rm where}\;\; 
		h_{\bf a}=
		\textstyle\left({\bf 1}_{\m F_p}+\sum_{\chi\in B}{\bf a}_\chi\chi\right)h_0\, .
	\end{equation*}
	These coordinates ${\bf a}=({\bf a}_{\chi})\in \m C^{B}$ are also a linear coordinate system
	for the tangent space of $\m P(V)$ at the point $[h_0]$, thanks to the formula
	\begin{equation*}
		\label{eqntanpv2}
		{\bf a}\mapsto v_{\bf a}:=\frac{d}{d\eps}[h_{\eps{\bf a}}]|_{\eps=0}
		\;\in\; T_{[h_0]}\m P(V) .
	\end{equation*}
	
	$\star$ {\bf We describe the tangent space to the  torus $T$}.
	Since the characters $\chi$ vanish at $0$,  
	the  linear condition defining the tangent space 
	of $T$ at the point $[h_0]$ is
	\begin{equation*}
		\label{eqntanf02}
		\textstyle
		{\rm Re}(\sum_{\chi\in B}{\bf a}_\chi \chi) \;=\; 0.
	\end{equation*}
	By the linear independance of the characters, this 
	can be rewritten as
	\begin{equation}
		\label{eqntanf04}
		\ol{{\bf a}_{\ol \chi}}=-{\bf a}_\chi
		\; \;\mbox{\rm for all $\chi\in B$}.
	\end{equation}
	
	$\star$ {\bf We describe the tangent space to the torus $F^{-1} T$}.
	Using \eqref{eqnchieig}, one computes in our coordinate system
	\begin{eqnarray*}
		h_{\bf a}
		&=&\textstyle
		h_0+\sum_{\chi\in B}
		(\cos \theta\, {\bf a}_\chi+
		i\, \sin\theta\,{\bf a}_{\chi_{_0}\chi})
		\chi\; ,\\
		\widehat{{h}}_{\bf a}&=&
		\textstyle
		\left({\bf 1}_{\m F_p}+{\bf b}_0\,\de_0+\sum_{\chi\in B}{\bf b}_\chi \ol\chi\right)\widehat{h}_0
		\;\;\;\;{\rm where}\\
		{\bf b}_0
		&=&  \tfrac{p-1}{\eps_p\sqrt{p}}\,
		(\cos \theta\, {\bf a}_{_{\bf 1}}+
		i\, \sin\theta\,{\bf a}_{\chi_{_0}})\, e^{i\theta_0},\\
		\textstyle
		\sum_{\chi\in B}{\bf b}_\chi \ol\chi
		\!&\!=\!&\!
		\textstyle
		\sum\limits_{\chi\in B}\tfrac{G(\chi)}{\eps_p\sqrt{p}}
		(\cos \theta\, {\bf a}_\chi\!+\!
		i\, \sin\theta\,{\bf a}_{\chi_{_0}\chi})
		(\cos \theta\, \ol\chi\! -\!
		i\, \sin\theta\,\chi_{_0}\ol\chi),
	\end{eqnarray*}
	where $G(\chi)$ is the Gauss sum \eqref{eqngausum}.
	The extra term ${\bf b}_0\,\de_0$ comes from the extra term in 
	\eqref{eqnchieig} when $\chi={\bf 1}_{\m F_p^*}$ is the trivial Dirichlet character.
	The coefficients ${\bf b}_\chi$ are then given by 
	\begin{equation}
		\label{eqnbchi}
		{\bf b}_\chi
		=
		\tfrac{G(\chi)\cos\theta^2+G(\chi_{_0}\chi)\sin\theta^2}{\eps_p\sqrt{p}}\,{\bf a}_\chi
		+
		i
		\tfrac{(G(\chi)-G(\chi_{_0}\chi))\cos\theta\,\sin\theta}{\eps_p \sqrt{p}}\,{\bf a}_{\chi_{_0}\chi}
	\end{equation} 
	In particular, for the trivial character, one has
	\begin{equation*}
		\label{eqnbchi1}
		{\bf b}_{\bf 1}
		\;:=\;
		\tfrac{-\cos\theta^2+\eps_p\sqrt{p}\,\sin\theta^2}{\eps_p\sqrt{p}}\,{\bf a}_{\bf 1}
		-
		i\, 
		\tfrac{(1+\eps_p\sqrt{p}\,)\,\cos\theta\,\sin\theta}{\eps_p \sqrt{p}}\,{\bf a}_{\chi_{_0}}
	\end{equation*}
	
	As above, the  linear condition defining the tangent space 
	of $F^{-1}T$ at the point $[h_0]$ is
	\begin{eqnarray}
		\label{eqntanf05}
		\ol{{\bf b}_{\ol \chi}}+{\bf b}_\chi
		&=& 0
		\; ,\;\;\mbox{\rm for all $\chi\in B$ non trivial, and}\\
		\nonumber\label{eqntanf06}
		{\rm Re}({\bf b}_{0}-{\bf b}_{\bf 1})
		&=& 0
		\; ,\;\;\mbox{\rm for the trivial character $\chi={\bf 1}_{\m F_p^*}$}.
	\end{eqnarray}
	
	$\star$ {\bf We check the transversality of these tangent spaces}.
	We want to prove that if ${\bf a}=({\bf a}_\chi)_{\chi\in B}$ belongs to both tangent spaces, then ${\bf a}=0$. 
	By \eqref{eqntanf04} one has
	$\ol{{\bf a}_{\ol \chi}}=-{\bf a}_\chi$
	for all $\chi\in B$. 
	One computes  $
	$ 
	\begin{equation}
		\label{eqnbchia}
		\hspace*{-1em}\chi(\!-\!1\!)\ol{{\bf b}_{\ol\chi}} =\!-\tfrac{G(\!\chi)\cos\!\theta^2+\chi_0(\!-\!1)G(\!\chi_{_0}\!\chi)\sin\!\theta^2}{\eps_p\sqrt{p}}{{\bf a}_{\chi}}
		\!+\! i
		\tfrac{(G(\!\chi)-\chi_0(\!-\!1)G(\!\chi_{_0}\!\chi))\cos\!\theta\,\sin\!\theta}{\eps_p\sqrt{p}}{{\bf a}_{\chi_{_0}\chi}}.\!\!\!
	\end{equation} 
	In order to exploit \eqref{eqntanf05},
	we distinguish two cases.
	\vs 
	
	{\bf First case: when $p\equiv 1$ mod $4$}.\\ 
	One has $\eps_p=1$,\; $\chi_0(-1)=1$,\; $\theta_0=0$ and $\cos\theta =\frac{1}{\sqrt{p}+1}$.

	\noindent 
	$\bullet$ If $\chi$ is odd. One combines the equations \eqref{eqnbchi}, \eqref{eqntanf05} and \eqref{eqnbchia}. One gets 
	$$
	(G(\chi)\cos\theta^2+G(\chi_0\chi)\sin\theta^2)\,{\bf a}_\chi=0.
	$$
	Hence since $|G(\chi)|=|G(\chi_0\chi)|$,
	one concludes  ${\bf a}_\chi=0$.
	
	\noindent 
	$\bullet$ If $\chi$ is even and $\chi\neq \chi_0$. One combines the equations \eqref{eqnbchi}, \eqref{eqntanf05} and \eqref{eqnbchia} for the character $\chi_0\chi$. One gets 
	\begin{equation}
		\label{eqngchgch}
		(G(\chi)-G(\chi_0\chi))\, {\bf a}_\chi=0. 
	\end{equation}
	Hence by Lemma \ref{lemgauquo}, one concludes ${\bf a}_\chi=0$.
	
	\noindent 
	$\bullet$ If $\chi=\chi_0$. Instead of  \eqref{eqngchgch}, one writes
	$0={\rm Re}({\bf b}_0-{\bf b}_{\bf 1})=i\sqrt{p}\,\sin\theta\, {\bf a}_{\chi_0}$. One concludes
	${\bf a}_{\chi_0}=0$.
	\vs 
	
	{\bf Second case: when $p\equiv 3$ mod $4$}.\\ 
	One has 
	$\eps_p=i$,\; $\chi_0(-1)=-1$,\; $\theta_0=\theta$
	and $\tan\theta =\sqrt{p}$. 
	
	\noindent 
	$\bullet$ If $\chi$ is odd and $\chi\neq \chi_0$. One combines the equations \eqref{eqnbchi}, \eqref{eqntanf05} and \eqref{eqnbchia}. One gets 
	\begin{equation}
		\label{eqnachach1}
		G(\chi)\cos\theta \, {\bf a}_\chi- i\, G(\chi_0\chi)\sin\theta
		\,{{\bf a}_{\chi_{_0}\chi}}
		\; =\; 0
	\end{equation} 
	Applying the same calculation to the even character $\chi_0\chi$, one gets
	\begin{equation}
		\label{eqnachach2}
		i \,G(\chi_0\chi)\cos\theta \, {\bf a}_\chi + G(\chi)\sin\theta\,{{\bf a}_{\chi_{_0}\chi}}
		\; =\; 0
	\end{equation}  
	The determinant of this linar system is 
	$D:=(G(\chi)^2-G(\chi_0\chi)^2)\cos\theta\sin\theta$.
	By Lemma \ref{lemgauquo}, this determinant is non zero. 
	Therefore one has   ${\bf a}_\chi=0$.
	
	\noindent 
	$\bullet$ If $\chi=\chi_0$. 
	Equation \eqref{eqnachach1} is still valid and gives
	${\bf a}_{\bf 1}+{\bf a}_{\chi_{_0}}=0$.
Instead of \eqref{eqnachach2}, one writes 
$0={\rm Re}({\bf b}_0-{\bf b}_{\bf 1})=-i\sqrt{p}\,{\bf a}_{\bf 1}$. Combining these two equalities, one concludes
${\bf a}_{\chi_{_0}}=0$.

\noindent 
$\bullet$ If $\chi$ is even. One applies the previous discussion to the odd character $\chi_0\chi$ and one also gets  ${\bf a}_\chi=0$.
\end{proof}

We have used the following

\bl
\label{lemgauquo}
Let $p$ be a prime number, let $\chi_0$ be the Dirichlet character on $\m F_p$ of order $2$ and let $\chi$ be another Dirichlet character on $\m F_p$.
Then the ratio of Gauss sums
$R(\chi):=\displaystyle\frac{G(\chi_0\chi)}{G(\chi)}$
is not a root of unity.
\el

\begin{proof} When $\chi=\chi_0$ or $\chi$ is trivial the ratio 
does not have absolute value $1$. Otherwise,
one has the equality 
$R(\chi)=\displaystyle\frac{G(\chi_0)}{J(\chi,\chi_0)}$ where $J$ is the Jacobi sum \eqref{eqnjacsum}. 
By Proposition \ref{projacsum},  the ratio 
$\displaystyle \frac{\ol{R(\chi)}}{R(\chi)}=\pm \frac{J(\chi,\chi_0)}{J(\ol{\chi},\chi_0)}$
is not a root of unity. Therefore the ratio $R(\chi)$ is not a root of unity either.
\end{proof}

This ends the proof of Theorem \ref{thmexibiu}.

{\small
\bibliography{fourier}
	}
	\vs 
	
	{\small
		\noindent
		Y. \textsc{Benoist}: CNRS, 
		Universit\'e Paris-Saclay,\hfill
		\texttt{yves.benoist@u-psud.fr}}
	
\end{document}